\documentclass[preprint]{elsarticle}
\usepackage{amsmath,amsfonts,amssymb,amsthm,latexsym}
\usepackage{tikz,verbatim}
\usepackage[dvipdf]{epsfig}
\usepackage{color}
\usepackage{graphicx,epsfig}
\usepackage{caption}
\usepackage{subfigure}
\let\today\relax
\makeatletter
\def\ps@pprintTitle{%
\let\@oddhead\@empty
\let\@evenhead\@empty
\def\@oddfoot{}%
\let\@evenfoot\@oddfoot}
\makeatother

\usepackage[colorlinks=true,citecolor=blue]{hyperref}
\newcommand{\textcite}[1]{\cite{#1}}

\newcommand{\be}{\begin{equation}}
\newcommand{\ee}{\end{equation}}

\newcommand{\dd}{\mathrm{d}}

\renewcommand{\vec}[1]{\boldsymbol{#1}}

\definecolor{remarktcolor}{HTML}{375502}
\definecolor{remarkcolor}{RGB}{193,42,132} 
\newcommand{\eps}{\varepsilon}
\allowdisplaybreaks

\begin{document}


\begin{frontmatter}
\title{Sharp-interface problem of the Ohta-Kawasaki model for symmetric diblock copolymers}
\author[address1]{Amlan K. Barua\corref{mycorrespondingauthor}}
\cortext[mycorrespondingauthor]{Corresponding author}
\ead{abarua@iitdh.ac.in}

\author[address2]{Ray Chew}

\author[address3]{Shuwang Li}

\author[address4]{John Lowengrub}

\author[address5]{Andreas M\"unch}

\author[address6]{Barbara Wagner}

\address[address1]{Department of Mathematics, IIT Dharwad, Karnataka, 580011, India}
\address[address2]{FB Mathematik \& Informatik, Freie Universit{\"a}t Berlin, Arnimallee 6, 14195 Berlin, Germany}
\address[address3]{Department of Applied Mathematics, Illinois Institute of Technology, Chicago, IL 60616, USA}
\address[address4]{Department of Mathematics, University of California Irvine, Irvine, California 92697, USA}
\address[address5]{Mathematical Institute, University of Oxford, Andrew Wiles Building, Woodstock Road, Oxford OX2 6GG, United Kingdom}
\address[address6]{Weierstrass Institute, Mohrenstrasse 39, 10117 Berlin, Germany}

\begin{abstract}
The Ohta-Kawasaki model for diblock-copolymers is well known to the scientific community of diffuse-interface methods. To accurately capture the long-time evolution of the moving interfaces, we present a derivation of the corresponding sharp-interface limit using matched asymptotic expansions, and show that the limiting process leads to a Hele-Shaw type moving interface problem. The numerical treatment of the sharp-interface limit is more complicated due to the stiffness of the equations. 
To address this problem, we present a boundary integral formulation corresponding to a sharp interface limit of the Ohta-Kawasaki model. 
Starting with the governing equations defined on separate phase domains, we develop boundary integral equations valid for multi-connected domains in a 2D plane. For numerical simplicity we assume our problem is driven by a uniform Dirichlet condition on a circular far-field boundary. The integral formulation of the problem involves both double- and single-layer potentials due to the modified boundary condition. In particular, our formulation allows one to compute the nonlinear dynamics of a non-equilibrium system and pattern formation of an equilibrating system. Numerical tests on an evolving slightly perturbed circular interface (separating the two phases) are in excellent agreement with the linear analysis, demonstrating that the method is stable, efficient and spectrally accurate in space. 
\end{abstract}


\end{frontmatter}


\section{Introduction}
The Ohta-Kawasaki (OK) model \cite{ohta1986equilibrium} was originally derived by Takao Ohta and Kyozi Kawasaki to investigate mesoscopic phase separation in block copolymers. The phase separation in copolymeric substances results in the formation of two distinct regions, each rich in a particular ingredient. Domains of various shape may emerge in the system under various ratios of molecular weight of the two species. It is necessary to investigate such systems as the resulting properties  are different from those observed in multiphase systems of single monomer types. The model has garnered strong interest since its emergence and has been connected to areas beyond which it was originally proposed. Examples of applications include problems in condensed matter physics and biological systems \cite{abels2018model}.

In the original work of Ohta and Kawasaki \cite{ohta1986equilibrium}, an energy functional was proposed to investigate the phenomenon of phase separation where both attractive (short-range) and repulsive (long-range) forces play their part in determining the configurations. The evolution equation corresponding to the functional and its steady version was first mentioned in \cite{nishiura1995}, where a connection was made between Hele-Shaw (HS) flow equations and the time-dependent OK problem. 
In this paper, we present a formal derivation of the corresponding sharp-interface limit using matched asymptotic expansions, and show that the limiting process leads to an HS-type moving interface problem.
This allows us to recast the long-time evolution of the OK problem as a modified HS problem and focus our attention to the latter to obtain insight into the original pattern formation problem. The analytical solutions are ruled out owing to the complicated geometry and we investigate the problem mainly using numerical approaches.   

The boundary integral method is a preferred choice as a numerical method for HS-type problems because it entails dimension reduction, i.e., the problem defined on a domain becomes a problem defined on the domain boundary. However, the equations of dynamics constitutes stiff equations due to the surface tension acting at the fluid-fluid interface, and without the special numerical techniques described in \cite{HLS}, it is practically impossible to perform long-time numerical simulations. Several references have used this technique with great success and we refer the interested reader to \cite{zhao2020pattern, zhao2018computation, ShuwangPRL, ShuwangPhysicaD, hou2001boundary}. 
We also note that our equations differ from the traditional HS equations \cite{ShuwangJCP} in a few subtle ways. In the  original HS model, the far-field boundary condition is of Neumann type which very naturally corresponds to injection/removal of the fluid. Our problem, on the other hand, is driven by a Dirichlet type boundary condition in the far-field. This renders the constraint on the integral of velocity to be different in our case. We also note that the far-field boundary is at a finite distance from the origin in our case while in the classical HS problems, the radius of the far-field boundary is infinite. 

The main contribution of this paper can be summarized as follows: starting with a rescaled formulation of the OK equation, we present a matched asymptotic analysis in the long-time limit that governs the dynamics of the emerging interfaces and this leads to modified HS equations of the OK model. We then prescribe a transformation that converts the HS equations from the Poisson equation to the Laplace equation and transform the interfacial and far-field boundary conditions accordingly. The equations are then investigated using a linear analysis. We prescribe a
boundary integral formulation for the Laplace equation using free-space Green's function and we investigate the boundary integral equations numerically as the analytical solutions are known in very limited cases. The numerical methods allow us to investigate the steady-state configuration for various patterns hitherto not explored in detail. Throughout our computation, we demonstrate high accuracy which is a trademark of boundary integral computations. Nonlinear computations indicate that the interface morphologies depend strongly on the mass flux into the system before the system reaching equilibrium. Simulations of multiple equilibrating interfaces show complicated interactions between phase domains including interface alignment and coarsening. 

This paper is organized as follows: In Section~\ref{gov_eq}, we give a formulation for the boundary value problem of the OK equation in a rescaled form that is suitable for the asymptotic analysis using matched asymptotic expansions, which is carried out in Section~\ref{sec:sil}. In Section~\ref{analytical}, the analytical solutions of the problem are discussed. Numerical methods on the boundary integral equations, the spatial discretization of the integral equations using spectrally accurate quadrature rules, the dynamical equations, and the small-scale decomposition are discussed in Section~\ref{num_methods}. The interface is updated based on these methods. Finally, we present results of numerical simulations in Section~\ref{results} and summarize our findings in Section~\ref{concl}.

\section{Formulation of the Ohta-Kawasaki phase-field model} \label{gov_eq}
In the framework of density functional theory, the OK problem in its dimensionless form it is \cite{ohta1986equilibrium, nishiura1995}
\begin{equation}\label{eqn:dft_okeqn_functional}
{\cal F}_{\textrm{OK}}[\phi]=\int_\Omega \frac12 (\nabla \phi)^2
+F(\phi)-F(\phi_-)+\frac\alpha 2 \psi (\phi-\bar\phi) \,\, dxdy.
\end{equation}
In a domain $\Omega$, $\phi(t,\vec{x})$ is the density difference, $\phi_A(\vec{x}) - \phi_B(\vec{x})$, at position $\vec{x} = (x,y)$ and at time $t$, where the overbar denotes the average of a quantity, e.g.\
\begin{equation}
\bar\phi\equiv\frac1{|\Omega|}\int_\Omega \phi \,\, dxdy\, .
\end{equation}
$\psi$ is given by the solution of the Poisson problem,
\begin{subequations}\label{psidef}
\begin{alignat}{2}
-\Delta \psi&=\phi-\bar\phi &&\qquad \text{on } \Omega,\\
\frac{\partial \psi}{\partial n_{\partial\Omega}} &= 0 &&\qquad \text{on } \partial \Omega,\\
\bar\psi&=0,
\end{alignat}
\end{subequations}
where the last condition is introduced to enforce the uniqueness of $\psi$.
Here, we use for the double-well free energy $F$ the form 
\begin{equation}
F(\phi)=\frac14\phi^4-\frac12\phi^2
\end{equation}
which has two minima at $\phi_\pm=\pm 1$.
%
The chemical potential $\mu$ is obtained by the first variation of the functional ${\cal F}_{\text{OK}}$
\begin{subequations}\label{OK2}
\begin{align}
\mu &= - \Delta \phi+(\phi^3-\phi)-\alpha\psi,
\intertext{which yields the flux}
j &= -\nabla \mu.
\intertext{The system is closed via mass conservation}
\frac{\partial \phi}{\partial t}&=-\nabla\cdot j
\label{OK2mc}
\end{align}
together with boundary and initial conditions 
\begin{align}
j\cdot n_{\partial\Omega}&=0, \qquad \frac{\partial \phi}{\partial n_{\partial\Omega}}=0
\qquad \text{on } \partial \Omega, \\
\phi(x,0) &= \phi_{\text{init}}(x).
\end{align}
\end{subequations}
Derivations of the Ohta-Kawasaki phase-field model using the gradient flow approach can be found in, e.g., \cite{zhang2006periodic, parsons2012numerical, le2010convergence}.

\section{The sharp-interface limit}\label{sec:sil}
For diblock copolymers, the long-time interface formation during phase separation that  sets the small-scale related to the interface width is directly connected to the parameter $\alpha$, via $\varepsilon=\alpha^{1/3}$ \cite{ohnishi1999analytical, choksi_2001}.
It is thus convenient to rescale the Ohta-Kawasaki model to this regime via 
 $x=\alpha^{-1/3}\tilde x$, $\psi=\alpha^{-2/3}\tilde \psi$,
$\mu=\alpha^{1/3}\tilde \mu$, $\tau=\alpha t$,
and ${\cal{\tilde F}}_{\textrm{OK}}=\varepsilon{\cal{ F}}_{\textrm{OK}}$. After dropping the tildes, the rescaled free energy can be written as 
\begin{equation}\label{FOK1re}
{\cal{ F}}_{\textrm{OK}}[\phi]=\int_\Omega \frac12 \eps (\nabla \phi)^2
+\eps^{-1}\left(F(\phi)-F(\phi_-)\right)+\frac1 2 \psi (\phi-\bar\phi),
\end{equation}
and thus the corresponding phase-field model 
\begin{subequations}\label{OK1re}
\begin{align}
\frac{\partial \phi}{\partial \tau}&=\Delta\mu,\\
\mu&=-\eps \Delta \phi+\eps^{-1}(\phi^3-\phi)-\psi,\\
-\Delta \psi&=\phi-\bar\phi, \\
\frac{\partial \phi}{\partial n_{\partial\Omega}} &= 0, 
\quad
\frac{\partial \psi}{\partial n_{\partial\Omega}} = 0,
\quad 
\frac{\partial \mu}{\partial n_{\partial\Omega}} = 0  \qquad \text{on } \partial \Omega,\\
\phi(x,0) &= \phi_0(x).
\end{align}
\end{subequations}
Due to the small parameter $\varepsilon$ multiplying the Laplace operator in the chemical potential, the problem is singularly perturbed as $\varepsilon\to 0$. While such problems have been considered before with different methods \cite{henry_2001, le_2008, nishiura1995}, we investigate this ``outer'' problem through matched asymptotic expansions, where asymptotic approximations for the outer problem are matched to approximations of a corresponding ``inner'' problem in the neighborhood of the sharp interface. Our investigation follows a similar method applied by \cite{pego_1989} for the Cahn-Hilliard equations. We assume $\phi(\tau,\vec{x})$, $\mu(\tau,\vec{x})$, and $\psi(\tau,\vec{x})$ have the asymptotic expansions,
$\phi  = \phi_0 + \varepsilon \phi_1 + \varepsilon^2 \phi_2 + 
\mathcal{O}(\varepsilon^3)$,
$\mu  = \mu_0 + \varepsilon \mu_1 + \varepsilon^2 \mu_2 + 
\mathcal{O}(\varepsilon^3)$, and 
$\psi  = \psi_0 + \varepsilon \psi_1 + \varepsilon^2 \psi_2 + 
\mathcal{O}(\varepsilon^3)$.
Substitution into \eqref{OK1re} yields the asymptotic problems for $\phi_i$ up to order $\varepsilon^2$,
\begin{equation}
\mathcal{O}\left(\varepsilon^{0}\right): \partial_\tau \phi_0 = \Delta \mu_0,\qquad
\mathcal{O}\left(\varepsilon^{1}\right): \partial_\tau \phi_1 = \Delta \mu_1,\qquad
\mathcal{O}\left(\varepsilon^{2}\right): \partial_\tau \phi_2 = \Delta \mu_2.
\end{equation}
Similarly for $\mu_i$,
\begin{subequations}
\begin{alignat}{2}
&\mathcal{O}\left(\varepsilon^{-1}\right):
\qquad 0 &&= F^\prime \left( \phi_0 \right),
\label{eqn:sil_outerChem1}\\
&\mathcal{O}\left(\varepsilon^{0}\right):
\qquad \mu_0 &&= F^{\prime\prime} \left( \phi_0 \right) \phi_1  + \psi_0,
\label{eqn:sil_outerChem2}\\
&\mathcal{O}\left(\varepsilon^{1}\right):
\qquad \mu_1 &&= F^{\prime\prime} \left( \phi_0 \right) \phi_2 + \frac{1}{2} F^{\prime\prime\prime}\left(\phi_0\right)\phi_1^2 - \Delta \phi_0 + \psi_1,
\label{eqn:sil_outerChem3}
\end{alignat}
\end{subequations}
and $\psi_i$,
\begin{equation}
\mathcal{O}\left(\varepsilon^{0}\right): - \Delta \psi_0 = \phi_0 - \bar{\phi},\quad
\mathcal{O}\left(\varepsilon^{1}\right): - \Delta \psi_1 = \phi_1,\quad
\mathcal{O}\left(\varepsilon^{2}\right): - \Delta \psi_2 = \phi_2.
\label{eqn:sil_outerNL3}
\end{equation}
On the fixed boundary $\partial\Omega$, the rescaled boundary conditions are 
\begin{equation*}
\frac{\partial\phi_i}{\partial n_{\partial\Omega}} = 0, \qquad \frac{\partial \mu_i}{\partial n_{\partial\Omega}} = 0, \qquad \frac{\partial\psi_i}{\partial n_{\partial\Omega}} = 0, \quad \text{on} \quad \partial\Omega \quad \text{for} \quad i=0,1,2,\dots
\label{eqn:sil_bcsOuter}
\end{equation*}

To derive the inner problems, it is convenient to introduce a parametrization $\vec{r}(\tau,s) = (r_1(\tau,s),r_2(\tau,s))$ of the free interface $\Gamma$ via the arc length $s$, and $\vec{\nu}(\tau,s)$, the normal inward-pointing vector along the free boundary, so that any point in the thin $\varepsilon$-region around $\Gamma$ can be expressed by  
\begin{equation*}
\vec{x}(\tau,s,z) = \vec{r}(\tau,s) + \varepsilon z \vec{\nu}(\tau,s), 
\end{equation*}
where $\varepsilon z$ is the distance along the inward normal direction $\vec{\nu}(\tau,s)$ from the sharp interface $\Gamma$, given by 
\begin{equation*}
\vec{\nu}(\tau,s) = \left( -\partial_s r_2, \partial_s r_1 \right),\quad
\vec{t}(\tau,s) = \left( \partial_s r_1, \partial_s r_2 \right).
\label{eqn:tvector}
\end{equation*}
The relation between the derivatives of a quantity $\tilde{v}(\tau,s,z)$ defined in inner coordinates and the derivatives in outer coordinates $v(\tau,\vec{x})$ can be expressed as a product of matrices, see \ref{sec:sil_derivation} and \cite{dreyer_wagner_2005}.

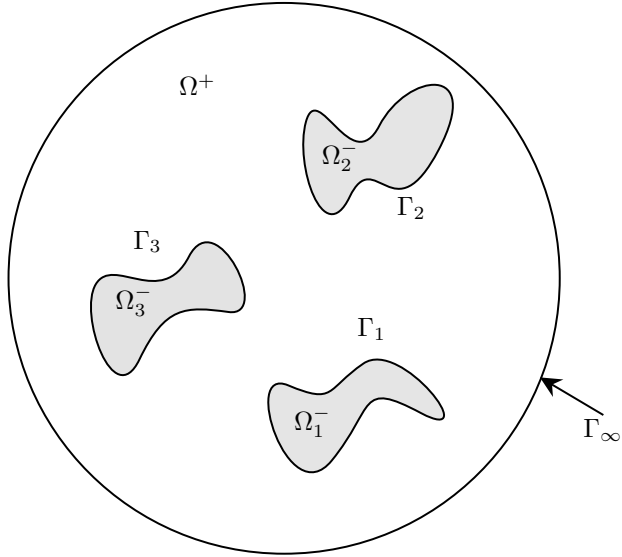
\begin{figure}
\centering
\tikzset{every picture/.style={line width=0.75pt}} 
\begin{tikzpicture}[x=0.75pt,y=0.75pt,yscale=-1,xscale=1]

\draw   (133.5,215) .. controls (133.5,138.23) and (195.73,76) .. (272.5,76) .. controls (349.27,76) and (411.5,138.23) .. (411.5,215) .. controls (411.5,291.77) and (349.27,354) .. (272.5,354) .. controls (195.73,354) and (133.5,291.77) .. (133.5,215) -- cycle ;
\draw  [fill={rgb, 255:red, 211; green, 211; blue, 211 }  ,fill opacity=0.62 ] (224.5,204) .. controls (238.5,177) and (266.5,235) .. (244.5,232) .. controls (222.5,229) and (212.5,229) .. (199.5,256) .. controls (186.5,283) and (168.5,233) .. (177.5,218) .. controls (186.5,203) and (210.5,231) .. (224.5,204) -- cycle ;
\draw  [fill={rgb, 255:red, 211; green, 211; blue, 211 }  ,fill opacity=0.58 ] (320.5,138) .. controls (334.5,111) and (373.5,106) .. (350.5,150) .. controls (327.5,194) and (318.5,147) .. (305.5,174) .. controls (292.5,201) and (277.5,158) .. (283.5,136) .. controls (289.5,114) and (306.5,165) .. (320.5,138) -- cycle ;
\draw  [fill={rgb, 255:red, 211; green, 211; blue, 211 }  ,fill opacity=0.57 ] (312.5,258) .. controls (331.5,245) and (373.5,300) .. (341.5,283) .. controls (309.5,266) and (318.5,280) .. (296.5,307) .. controls (274.5,334) and (250.5,259) .. (274.5,269) .. controls (298.5,279) and (293.5,271) .. (312.5,258) -- cycle ;
\draw    (433.5,284) -- (404.08,266.53) ;
\draw [shift={(401.5,265)}, rotate = 390.7] [fill={rgb, 255:red, 0; green, 0; blue, 0 }  ][line width=0.08]  [draw opacity=0] (10.72,-5.15) -- (0,0) -- (10.72,5.15) -- (7.12,0) -- cycle    ;

\draw (276,278.4) node [anchor=north west][inner sep=0.75pt]    {$\Omega ^{-}_{1}$};
\draw (290,144.4) node [anchor=north west][inner sep=0.75pt]    {$\Omega ^{-}_{2}$};
\draw (186,217.4) node [anchor=north west][inner sep=0.75pt]    {$\Omega ^{-}_{3}$};
\draw (218,110.4) node [anchor=north west][inner sep=0.75pt]    {$\Omega ^{+}$};
\draw (308,233.4) node [anchor=north west][inner sep=0.75pt]    {$\Gamma _{1}$};
\draw (195,189.4) node [anchor=north west][inner sep=0.75pt]    {$\Gamma _{3}$};
\draw (328,172.4) node [anchor=north west][inner sep=0.75pt]    {$\Gamma _{2}$};
\draw (422,284.4) node [anchor=north west][inner sep=0.75pt]    {$\Gamma _{\infty} $};
\end{tikzpicture}  
\caption{A schematic diagram of Ohta-Kawasaki problem. The interior domain $\Omega^-$ is the disjoint union of three connected and bounded regions $\Omega^{-}_1, \Omega^{-}_2$ and $\Omega^{-}_3.$ The boundary  of $\Omega^-$ consists of $\Gamma=\partial \Omega^-_1 \cup \partial \Omega^-_2 \cup \partial \Omega^-_3.$ The outer region $\Omega^+$ is bounded and surrounds $\Omega^-$.}
\label{fig_scheme}
\end{figure}

Similar to the outer problem, we assume that the inner asymptotic expansions for $\tilde{\phi}(\tau, s, z)$, $\tilde{\mu}(\tau, s,z )$, and $\tilde{\psi}(\tau, s,z )$ are given by
$\tilde{\phi} = \tilde{u}_0 + \varepsilon \tilde{u}_1 + \varepsilon^2 \tilde{u}_2 + \mathcal{O}(\varepsilon^3)$, 
$\tilde{\mu} = \tilde{\mu}_0 + \varepsilon \tilde{\mu}_1 + \varepsilon^2 \tilde{\mu}_2 + \mathcal{O}(\varepsilon^3)$, and
$\tilde{\psi} = \tilde{\psi}_0 + \varepsilon \tilde{\psi}_1 + \varepsilon^2 \tilde{\psi}_2 + \mathcal{O}(\varepsilon^3)$.
After application of the coordinate transformations to the governing equations, we obtain asymptotic subproblems for $\tilde\phi$, $\tilde\mu$ and $\tilde\psi$ for the inner region. These problems are solved and matched to the outer solutions. The details of the arguments, the matching conditions for the asymptotic analysis, are carried out in \ref{sec:sil_derivation}, resulting in the sharp-interface problem
\begin{subequations}\label{sil}
\begin{alignat}{2}
\phi_0 &= \pm 1,\\
 -\Delta \psi_0 &=\phi_0- \bar\phi && \qquad \text{ in } \Omega, \\
\Delta \mu_0 &= 0 && \qquad \text{ in } \Omega^\pm,\\
 \mu_0 &= \sigma\kappa-\psi_0 && \qquad \text{ on } \Gamma,  \\
 V &= \frac{1}{2} \left[ \frac{\partial \mu_0}{\partial n}\right] && \qquad \text{ on } \Gamma,\\
 \frac{\partial \mu_0}{\partial n_\infty} &= 0, \qquad
\frac{\partial \psi_0}{\partial n_\infty} = 0 && \qquad \text{ on } \partial \Omega,  
\end{alignat}
\end{subequations}
where $\sigma$ is the surface tension and  $\Omega = \Omega^+ \cup \Gamma \cup \Omega^-$ a domain, with $\Omega^+$ and the $\Omega^-$ the regions where $\phi_0=+1$ and $\phi_-=-1$, respectively, and $\Gamma$ is the interface between them. The normal to the latter pointing from $\Omega^+$ to $\Omega^-$ is called $n$. We will, more specifically, denote by $\Omega^+$ the exterior and $\Omega^-$ the interior domain. The boundary of $\Omega$ is denoted by $\partial \Omega$ and the jump of $\mu$ across the interface $\Gamma$ is given by
\[
\left[ \frac{\partial \mu_0}{\partial n}\right] = \frac{\partial \mu_0^+}{\partial n} - \frac{\partial \mu_0^-}{\partial n}.
\]
Finally, the value of $\sigma$ can be expressed as  
\begin{equation}
\sigma=\frac{1}{\phi_+-\phi_-}\int_{\phi_-}^{\phi_+} \sqrt{2 (F(\phi)-F(\phi_-))} \ \mathrm d\phi\, .
\end{equation}
For the derivation of the boundary integral formulation, it is convenient to reformulate the sharp-interface problem in terms of the variable
\begin{equation}
u:= \psi_0 + \mu_0. \label{bu}
\end{equation}
We consider a bounded domain $\Omega = \Omega^+ \cup \Gamma \cup \Omega^- \subset \mathbb{R}^2$ where $\Omega^+$, the outer domain, and $\Omega^-$, the inner domain, are open sets of $\mathbb{R}^2$ and $\Gamma$ is the moving interface separating the exterior domain $\Omega^+$ and the interior domain $\Omega^-$. The interior domain $\Omega^-$ is a disjoint union of finitely many open, connected components $\Omega_1^{-}, \Omega_2^{-}, \cdots, \Omega_
M^{-}$ and thus $\Gamma = \cup_{k=1}^M \partial \Omega_k^{-}.$  
The outer boundary of $\Omega$ is denoted by $\Gamma_{\infty}$. A schematic diagram of the problem is given in Fig.~(\ref{fig_scheme}).
The sharp-interface model is the following problem:
\begin{subequations}
\begin{alignat}{2}
-\Delta u &= 1-2\chi_{\Omega^-} && \qquad \text{ in } \Omega \backslash \Gamma, \label{Poisson1}\\
u &= \sigma\kappa && \qquad \text{ on } \Gamma, \label{bc1} \\
\frac{\partial u}{\partial n_\infty} &= 0 && \qquad \text{ on } \Gamma_{\infty}, \label{bc2} \\
V &= \frac{1}{2} \left[ \frac{\partial u}{\partial n}\right] && \qquad \text{ on } \Gamma, \label{vel}
\end{alignat}
\end{subequations}
where $u$ is an unknown function, $\chi_A$ is the characteristic function of the set $A$, $\kappa$ is the curvature of boundary $\Gamma$, $\sigma$ is the surface tension parameter, the operator $\dfrac{\partial}{\partial n}$ is the normal derivative where ${\bf n}$ denotes the normal directed from $\Omega^-$ to $\Omega^+$. While the function $u$ is continuous, the  derivative of $u$ suffers a jump across the interface $\Gamma$ and is given by $\left[ \dfrac{\partial u}{\partial n}\right] = \dfrac{\partial u^+}{\partial n} - \dfrac{\partial u^-}{\partial n} $, where $u^+$ and $u^-$ are the solutions of the OK problem in the  exterior and interior domains respectively.
The interface $\Gamma$ moves due to the velocity $V$.

To eliminate the source term in the field equation and recast the problem in terms of the Laplace equation, we introduce a new function $w$ defined as 
\begin{equation}
w = u + \frac{\left(1-2\chi_{\Omega^-}\right)}{4}\left|{\bf x}\right|^2 \label{w},
\end{equation}
where $\left|{\bf x}\right|^2 = x^2+y^2$.  Then the functions $u^+$ and $u^-$
are replaced by $w^+ = u^+ + \frac{1}{4}\left|{\bf x}\right|^2$ and $w^- = u^- - \frac{1}{4}\left|{\bf x}\right|^2$ in $\Omega^{+}$ and $\Omega^{-}$ respectively.  
The boundary condition Eq.~\eqref{bc1} on $\Gamma$ splits into conditions on $w^{-}$ and $w^{+}$ as follows:
\begin{align}
w^- &= \sigma\kappa - \frac{\left|{\bf x}\right|^2}{4},\\
w^+ &= \sigma\kappa + \frac{\left|{\bf x}\right|^2}{4}.
\end{align}
We also transform the far-field boundary condition Eq.~(\ref{bc2})  to
\begin{equation}
\label{newfar}
\frac{\partial w^+}{\partial n_\infty} = \frac{1}{2}{\bf x_\infty}\cdot {\bf n}_\infty,
\end{equation}
where ${\bf x}_\infty$ is a point on the outer boundary $\Gamma_{\infty}$ and ${\bf n}_\infty $ is the outward normal at ${\bf x}_\infty$.  The normal velocity of the interface $\Gamma$ separating the interior and the exterior domain becomes
\begin{equation}
V = \frac{1}{2} \left[ \frac{\partial w}{\partial n}\right] - \frac{1}{2}{\bf x}\cdot {\bf n} \label{newvel},
\end{equation}
where, as in Eq.~(\ref{vel}), $\left[ \dfrac{\partial w}{\partial n}\right] = \dfrac{\partial w^+}{\partial n} - \dfrac{\partial w^-}{\partial n}$.

\section{Analytical solution of original equations}\label{analytical}
It is not possible to find analytical solutions of the OK equations for arbitrary geometry and multiply connected regions. However, for simplified cases, like when $\Omega^- \cup \Gamma\cup \Omega^+$ is a circular domain centered at origin and $\Omega^-$  a circular domain of smaller radius and centered at zero, it is possible to find an analytical solution. In such a case \cite{Barua}, the solution inside $\Omega^{-}$ is obtained as 
\begin{equation}
u^{-} = \frac{1}{4}\left(x^2+y^2-R^2\right)+\frac{\sigma}{R}. \label{sol11} 
\end{equation}
Similarly, in the exterior domain, the solution of the boundary value problem of the Poisson equation in $\left(r,\theta\right)$ coordinates is given by
\begin{equation}
u^+(r) = -\frac{r^2}{4}+\left(\frac{R_{\infty}^2}{2} \right)\log r + \frac{\sigma}{R} + \frac{R^2}{4} - \frac{R_{\infty}^2}{2}\log R. \label{sol2} 
\end{equation}
In steady state, the interface between the two domains does not move ($V=0$)
and Eq.~\eqref{vel} requires the normal derivative of $u$ to be continuous. From Eq.~\eqref{sol11}
and \eqref{sol2}, we get
\begin{subequations}
\begin{align}
\left.\frac{\partial u}{\partial n}\right|_{R^-} &=\frac R2,\label{veleqA}\\
\left.\frac{\partial u}{\partial n}\right|_{R^+} &=-\frac R2+\frac{R_\infty^2}{2R}. \label{veleqB}
\end{align}
\end{subequations} 
Equating the two gives an additional relation between the radii of the interior 
and the total domain,
\begin{equation}
R_\infty=\sqrt2\, R,
\end{equation}
which simply states that the area of the interior and exterior domains are equal, as expected
for a symmetric diblock copolymer configuration in steady state.

The solution of the OK equations can be extended further via linear analysis on a domain $\Omega^{-}$ with the shape of a slightly perturbed circle of the form
\begin{equation}
 r\left(t, R,\theta\right) = R\left(t\right) + \delta\left(t\right) \cos k\theta, \quad 0\leq \theta <2\pi,
\end{equation}  
where $R$ is the radius of the circle and $\delta \cos k\theta$ is a small perturbation  with  $\dfrac{\delta\left(0\right)}{R\left(0\right)}\sim \mathcal{O}\left(\epsilon\right), \epsilon \ll 1$. Thus, by continuity of the problem, we expect $\dfrac{\delta\left(t\right)}{R\left(t\right)}\sim \mathcal{O}\left(\epsilon\right)$, at least for $t\leq T$, where $T>0$ is possibly a short period of time. In this case, it is easier to work with the transformed equations and 
we presume that the solution in polar coordinates is given by 
\begin{equation}
w^{\pm}\left(r,\theta\right)= w^{\pm}_0\left(r\right)+\delta w^{\pm}_1\left(r,\theta\right)+\mathcal{O}\left(\delta^2\right),
\end{equation}
where $w^{\pm}_0$ is the zeroth order solution and $w^{\pm}_1$ is the first order solution.
A straightforward computation yields the zeroth order solution as
\begin{subequations}
\begin{align}
 w^{-}_0 &= \frac{\sigma}{R} - \frac{R^2}{4},\label{sol0}\\
 w^{+}_0 &= \frac{R_{\infty}^2}{2}\log r +\frac{\sigma}{R}+\frac{R^2}{4} - \frac{R_{\infty}^2}{2}\log R. \label{sol1}
\end{align}
\end{subequations}
Next we compute the first order corrections and in this case, $w^-$ is of the form $A^{-}r^k \cos k\theta$ where 
\begin{equation}
A^{-} = \frac{\sigma\left(k^2-1\right)}{R^{k+2}}-\frac{1}{2R^{k-1}}.
\end{equation}
The function $w^+$ is of the form $\left[A^+r^k + \dfrac{B^+}{r^k}\right]\cos k\theta$ where
\begin{eqnarray}
 A^+ &=& \frac{R^k}{R^{2k}+R_{\infty}^{2k}}\left[ \frac{\sigma\left(k^2-1\right)}{R^2} + \frac{R}{2} -\frac{R^2_{\infty}}{2R}\right],\\
 B^+ &=& \frac{R^k R_{\infty}^{2k}}{R^{2k}+R_{\infty}^{2k}}\left[ \frac{\sigma\left(k^2-1\right)}{R^2} + \frac{R}{2}-\frac{R^2_{\infty}}{2R}\right].
\end{eqnarray}
Once the functions $w^-$ and $w^+$ are available up to first order, we may proceed to calculate the velocity of the interface as
\begin{equation}
   V\approx \dot{r} = \dot{R}+\dot{\delta}\cos k \theta \label{velrad}
\end{equation}
where the ``dot'' on the respective variables indicate derivative with respect to time. The expression on the right of Eq.~(\ref{velrad}) captures the interface velocity up to first order. We equate the right hand side of Eq.~(\ref{velrad}) to the right hand side of Eq.~(\ref{newvel}) and obtain
\begin{align}
\dot{R}&= R_{\infty}^2/4R-R/2, \label{R_t}\\
\dot{\delta}&= \left[-R_{\infty}^2/R^2 + k(t_2-t_3)/2 - k t_1/2-1/2\right]\delta. \label{d_t}
\end{align}
where 
\begin{align}
t_1 &= \sigma(k^2-1)/R^3 - 1/2,\\
t_2 &= p_1R^{2k-1}/(R^{2k} + R_{\infty}^{2k}),\\ 
t_3 &= p_1R_{\infty}^{2k}/(R(R^{2k} + R_{\infty}^{2k})),\\
p_1 &= \sigma(k^2-1)/R^2 + R/2 - R_{\infty}^2/\left(2R\right).
\end{align}
These solutions are used later on to validate our numerical methods.

\section{Numerical methods}\label{num_methods}
In this section, we describe the numerical methods including the derivation of the boundary integral equation, its solution, and methods to update the interface. The switch from differential equation to boundary integrals results in a dimension reduction as the original PDE problem should be solved over a domain while the integral equations only have to be solved on the boundary. 
\subsection*{Mathematical preliminaries}
We observe  that the interface $\Gamma$, on which we have to solve the integral equation, is a union of disjoint, smooth, and closed curves $\partial \Omega_k^-, k=1, \cdots, M$ where $\partial \Omega_k^-$ is the boundary of the region $\Omega_k^-.$ 
We assume that each interface $\partial \Omega^-_k$ is represented by
\begin{equation}
 \partial \Omega_k^- = \left\{ {\bf x}\left(\alpha,t\right)=\left(x\left(\alpha,t\right), y\left(\alpha,t\right)\right):0\leq \alpha<2\pi\right\},
\end{equation}
where the function ${\bf x}$ is analytic and $2\pi$-periodic in the parameter $\alpha$. The local tangent and the normal vectors to the interface are
\begin{equation}
{\bf s} = \left(x_{\alpha}, y_{\alpha}\right)/s_{\alpha} \quad \text{and} \quad {\bf n} = \left(y_{\alpha}, -x_{\alpha}\right)/s_{\alpha}
\end{equation}
respectively, where 
$x_\alpha$ and $y_\alpha$ are the derivatives w.r.t. to $\alpha$
and $s_{\alpha}=\sqrt{x_{\alpha}^2+y_{\alpha}^2}$ is the local variation of arc length. If we introduce the angle $\theta$ tangent to the interface, then we may write ${\bf n} = \left(\sin \theta, -\cos \theta\right)$ and the curvature $\kappa = \theta_{\alpha}/s_{\alpha} = \theta_s.$

\subsection*{Boundary integral formulation}\label{bie}
The introduction of the function $w$ in Eq.~(\ref{w}) allows us to transform the Poisson equation in the original problem to the Laplace equation. We further wish to recast the latter using boundary integral formulation.
Consider the free space Green's function $G({\bf x},{\bf x}^\prime)=\frac{1}{2\pi} \ln |{\bf x}-{\bf x}^\prime|$. We then write the solution $w^-$ to the interior problem as a combination of single layer and double layer potential, i.e., 
\begin{equation}
w^-\left({\bf x}\right) = \int_\Gamma \left\{\frac{\partial w^-\left({\bf x^\prime}\right)}{\partial n\left({\bf x^\prime}\right)} G\left({\bf x},{\bf x}^{\prime}\right) - w^-\left({\bf x^\prime}\right) \frac{\partial G}{\partial n\left({\bf x^\prime}\right)} \right\} ds^\prime,
\end{equation}
for ${\bf x} \in \Omega^-.$ As ${\bf x} \rightarrow {\bf x}^\prime \in \Gamma$, we have
\begin{equation}
\frac{1}{2}\left(\sigma\kappa - \frac{\left|{\bf x}\right|^2}{4}\right) = \int_\Gamma \left\{\frac{\partial w^-(\bf x^\prime)}{\partial n({\bf x^\prime})} G({\bf x},{\bf x}^\prime) - w^-\left({\bf x^\prime}\right) \frac{\partial G}{\partial n\left({\bf x^\prime}\right)} \right\} ds^\prime \label{insideongamma}.
\end{equation}
Similarly for the exterior problem,  
\begin{equation}
w^+({\bf x}) = \tilde w_\infty - \int_{\Gamma} \left\{\frac{\partial w^+({\bf x^\prime})}{\partial n({\bf x^\prime})} G({\bf x},{\bf x}^\prime) - w^+({\bf x^\prime}) \frac{\partial G }{\partial n({\bf x^\prime})} \right\} ds^{\prime},
\end{equation}
for ${\bf x} \in \Omega^+$, where $\tilde w_\infty$ is an unknown to be solved. 
As ${\bf x} \rightarrow {\bf x}^\prime \in \Gamma$, we have
\begin{equation}
\frac{1}{2}\left(\sigma\kappa + \frac{\left|{\bf x}\right|^2}{4}\right) =\tilde{w}_\infty - \int_\Gamma \left\{ 
\frac{\partial w^+\left({\bf x^\prime}\right)}{\partial n\left({\bf x^\prime}\right)} G({\bf x},{\bf x}^\prime) - w^+({\bf x^\prime}) \frac{\partial G}{\partial n(\bf x^\prime)} 
\right\} ds^\prime \label{outsideongamma}.
\end{equation}
Adding equations (\ref{insideongamma}) and (\ref{outsideongamma}) together, we have
 \begin{equation}
\sigma\kappa  =\tilde w_\infty - \int_\Gamma 2V G({\bf x},{\bf x}^\prime) ds^\prime - \int_\Gamma ({\bf x^{\prime}}\cdot {\bf n}^{\prime}) G({\bf x},{\bf x}^\prime) ds^\prime +\int_\Gamma \frac{|{\bf x}^{\prime}|^2}{2}  \frac{\partial G}{\partial n({\bf x^\prime})} ds^\prime \label{BIE}.
\end{equation}
Eq.~(\ref{BIE}) is the boundary integral equation that we solve numerically.  
An additional equation is needed to complete the problem.
To this end, we integrate $\Delta w^-=0$ in $\Omega^-$ and $\Delta w^+=0$ in $\Omega^+$, and we then use the divergence theorem to get $\int_\Gamma \frac{\partial w^-}{\partial n} ds=0$ and $\int_\Gamma \frac{\partial w^+}{\partial n} ds+ \int_{\partial \Omega} \frac{\partial w^+}{\partial n_\infty} ds=0$. Subtracting these two equations and using equation (\ref{newvel}), we get
\begin{equation} 
 J=\int_{\Gamma} V ds =  \frac{1}{2}A_{total} - A^- \label {bc3new},
\end{equation}
where $A_{total}$ is the total area enclosed by $\Gamma_{\infty}$ and $A^-$ is the area enclosed by $\Gamma$. 
We solve for $\tilde w_{\infty}$ and the normal velocity $V$ using equations (\ref{BIE}) and (\ref{bc3new}). The physical meaning of  $\tilde w_{\infty}$ in the integral equation is evident: It is the value of $w$ at $\Gamma_{\infty}$ corresponding to the flux given in the right hand side of Eq.~(\ref{bc3new}). Our formulation thus allows us to investigate the (unknown) Dirichlet condition at the far-field corresponding to a (known) Neumann condition.  
\subsection*{Solving the integral equations}
\noindent  The boundary integral equation (\ref{BIE}) in equal arc length parameter is given by  
\begin{align}
\tilde w_\infty &- \int_\Gamma 2V\left({\bf x}\left(\alpha^\prime\right) \right) G({\bf x}\left(\alpha\right),{\bf x}\left(\alpha^\prime)\right)s_{\alpha}\left(\alpha^\prime\right) d\alpha^\prime \nonumber\\
&=\sigma\kappa + \int_\Gamma \left({\bf x}\left(\alpha^\prime\right)\cdot {\bf n}\left(\alpha^\prime\right)\right) G({\bf x}\left(\alpha\right),{\bf x}\left(\alpha^\prime)\right) s_{\alpha}\left(\alpha^\prime\right) d\alpha^\prime \nonumber\\
&\quad -\int_\Gamma \frac{|{\bf x}\left(\alpha^\prime)\right)|^2}{2}  \frac{\partial G({\bf x}\left(\alpha\right),{\bf x}\left(\alpha^\prime)\right)}{\partial n({\bf x}\left(\alpha^\prime\right))} s_{\alpha}\left(\alpha^\prime\right) d\alpha^\prime.
\label{BIE3A}
\end{align}
This along with  Eq.~(\ref{bc3new})   should be solved to find the  velocity $V$ of the interface as well as $w_{\infty}$. We use the Nystr{\"o}m method to discretize the integral equations using highly accurate quadrature rules on the various integrals in Eq.~(\ref{BIE3A}). We discretize each of the curves using $N$ marker points using equal arc length parametrization $\alpha_j = jh$ where $h=2\pi /N$. We choose $N=2^n$ for some positive integer $n$.   
Next, we investigate the smoothness of the various integrals in Eq.~(\ref{BIE3A}). 
\subsubsection*{Double-layer potential}
The kernel $\frac{\partial G}{\partial n({\bf x^\prime})}$ of the integral 
$\int_\Gamma   \frac{\partial G}{\partial n({\bf x^\prime})} \frac{\left|{\bf x}\right|^2}{2} ds^\prime 
$
does not a have a singularity as $\dfrac{\partial \log\left| {\bf x}\left(\alpha\right)-{\bf x}\left(\alpha^{\prime}\right)\right|}{\partial n({\bf x}\left(\alpha^{\prime}\right))} = \frac{1}{2}\kappa\left(\alpha\right)+\mathcal{O}\left(\alpha-\alpha^{\prime}\right)$
with $\alpha^{\prime}\rightarrow \alpha.$
Thus, an application of trapezoidal or alternating point quadrature is enough to ensure spectral accuracy \cite{SI}. One may also apply the hybrid Gauss-trapezoid quadrature rules derived using the Euler-Maclaurin formula, as suggested in \cite{alpert1999hybrid}.
\subsubsection*{Single-layer potential}
The second integrals, both in the left and right hand side of Eq.~(\ref{BIE3A}), possess a logarithmic singularity and cannot be handled by trapezoidal rule as it is only second-order accurate. However, the integration can be performed by first splitting the log kernel  as
\begin{equation}
    \log \left|x\left(\alpha,t\right)-x\left(\alpha^{\prime},t\right)\right| = \log2\left|\sin\left(\frac{\alpha-\alpha^{\prime}}{2}\right) \right| +\log \frac{\left|x\left(\alpha,t\right)-x\left(\alpha^{\prime},t\right)\right|}{2\left|\sin\left(\frac{\alpha-\alpha^{\prime}}{2}\right) \right|},
\end{equation}
and then by applying the additive rule of integration.  
The kernel of the integration $\int_{0}^{2\pi} f\left(\alpha,\alpha^{\prime}\right)\log 2\left| \sin\left(\frac{\alpha-\alpha^{\prime}}{2}\right) \right|d\alpha^{\prime}$
 is  still singular at $\alpha=\alpha^{\prime}$, but the use of a Hilbert transform \cite{HLS} or quadrature referred in \cite{kress1995numerical} results in spectral accuracy. In this work we use the method suggested in \cite{HLS}. 
The kernel of second integration $\int_{0}^{2\pi} f\left(\alpha,\alpha^{\prime}\right) \log \dfrac{\left|x\left(\alpha,t\right)-x\left(\alpha^{\prime},t\right)\right|}{2\left|\sin\left(\frac{\alpha-\alpha^{\prime}}{2}\right) \right|}d\alpha^{\prime}$ has a removable singularity at $\alpha=\alpha^{\prime}$ and can be evaluated via alternating point quadrature rule.

The overall discretization of the integral equation gives rise to a dense system of linear equations comprising of $MN+1$ equations, where $M$ is the number of connected components of $\Omega^{-}$ and $N$ is the number of marker points on the boundary of each component. We have an additional unknown in the form of $w_{\infty}$.
 We solve this system using an iterative GMRES \cite{SS} technique. The GMRES requires only the (dense) matrix-vector multiplication routine and this is the most time consuming part of the iterative solver. Since our matrix is dense, the routine is completed by $\mathcal{O}\left(M^2N^2\right)$ operations. The cost of matrix-vector multiplication operation can be reduced by the application of a parallel matrix-vector multiplication. It can also be reduced to $\mathcal{O}\left(MN\log\left(MN\right)\right)$ by the use of fast summation algorithms \cite{Greengard1,feng2014parallel, Lindsay1}. We do not use any preconditioner in the solver.   
\subsection*{Evolution of domain interfaces}
The discretization of the integral equation gives rise to a stiff system of ODEs as the motion of the interface is curvature driven~\cite{HLS}. The time explicit methods result in a stability constraint $\Delta t \sim \mathcal{O}\left(\Delta s^3\right)$ where $\Delta s$ is the spatial resolution. Moreover, the Lagrangian marker points can come close to each other during the course of evolution.
To circumvent these problems, we  implement the small scale decomposition technique due to Hou et. al.~\cite{HLS}. This special temporal scheme reduces the stiffness requirement to $\Delta t \sim \mathcal{O}\left(\Delta s\right)$. The scheme also prevents two points from coming too close to each other by distributing the markers on the interface using equal arc length frame and then maintaining the same at all time by the addition of a tangential velocity $T$ at every step of calculation. 
\subsubsection*{Dynamics of the interface}
Once the velocity $V$ is obtained for each marker point, we do not 
update Eq.~(\ref{newvel}) directly. Instead, the dynamics of the problem is recast in terms of the lengths $L$ of the interfaces and the angle $\theta$ that the tangent to the marker point makes with the positive $x$-axis. First, we add a tangent velocity $T\left(\alpha,t \right)$ to the interface where $T\left(\alpha,t \right)$ is given by
\begin{equation}
T\left(\alpha,t\right) =    T\left(0,t\right) + \int_{0}^{\alpha} s^{\prime}_{\alpha} \kappa^{\prime} V d\alpha^{\prime}-\frac{\alpha}{2\pi}  \kappa^{\prime} V \, d\alpha^{\prime}. 
\end{equation}
After adding the tangential velocity, the motion of the interface is given by 
\begin{equation}
    \frac{d}{dt}{\bf x}\left(\alpha,t\right) = V\left(\alpha,t\right){\bf n}+T\left(\alpha,t\right){\bf s}.
\end{equation}
The addition of the tangential velocity does not change the shape of the interface; however, it is crucial for maintaining the equal arc length distribution of the marker points throughout the computation and prevents the clustering problem. Once the equal arc length distribution is taken care of, we pose the dynamics of the problem with the following two equations,
\begin{align}
L^i_t &= \int_0^{2\pi} \theta^i_{\alpha}V^i \left(\alpha,t\right) \, d\alpha, \label{lq}\\
\theta^{i}_{t} &= \frac{2\pi}{L^i}\left(-V^i_{\alpha} + T^i \theta^i_{\alpha}\right), \qquad i=1, \dots, M. \label{thetaq}
\end{align}
The subscripts $\alpha$ and $t$  denote derivatives with respect to these variables. We use an additional superscript $i$ to indicate the interface for which the equations are written. We obtain one equation for $L$ for each of the $M$ domains, while we get one equation for $\theta$ for every marker point on the boundaries of the domains. Thus, we must solve $M+MN$ ordinary differential equations in total. It should be noted that the interface can be fully recovered from $L$ and $\theta$ by integrating the relation
\begin{equation}
    {\bf x}_{\alpha}^i = \frac{L^i\left(t\right)}{2\pi}\left(\cos \theta^i\left(\alpha,t\right),\sin \theta^i\left(\alpha,t\right)\right).
\end{equation}
\subsubsection*{Small-scale decomposition and updating the interface}
The stiffness of the original problem propagates to Eq.~(\ref{thetaq}), while Eq.~(\ref{lq}) is non-stiff. The latter can be integrated explicitly, but the solution technique for the $\theta$-equation is far from trivial. This equation is solved using small-scale decomposition (SSD), an idea which has been successfully used in a number of problems in the domain of, e.g., HS flow, micro-structure evolution \cite{Shuwang2, leo2000microstructural}, vesicle wrinkling \cite{sohn2010dynamics}, and dynamics of an epitaxial island \cite{li2011boundary}. In problems driven by Laplace-Young boundary conditions, the critical factor in the numerical computation is the curvature of the interface. It introduces higher derivatives in the dynamical equations and results in severe stability constraints. For example, the analysis of the equations of motion reveals \cite{HLS} that, at small spatial scales,
$V\left(\alpha,t\right) \sim \frac{\sigma}{s^2_{\alpha}}\mathcal{H}\left[\theta_{\alpha\alpha}\right] $
where $\mathcal{H}\left[\theta_{\alpha\alpha}\right]$ denotes the periodic Hilbert transform of $\theta_{\alpha\alpha}$ and therefore Eq.~(\ref{thetaq}) becomes
\begin{equation}
    \theta_t = \frac{\sigma}{s^3_{\alpha}}\mathcal{H}\left[\theta_{\alpha\alpha\alpha}\right]+N\left(\alpha,t\right) \label{SSD},
\end{equation}
where the term $N\left(\alpha,t\right) = \left(V_s+\kappa T\right)-\frac{\sigma}{s^3_{\alpha}}\mathcal{H}\left[\theta_{\alpha\alpha\alpha}\right].$ In the last equation and the subsequent ones, we suppress $i$ in the superscript to keep our notation simple, but its presence should be understood. SSD reveals that the part $ \frac{\sigma}{s^3_{\alpha}}\mathcal{H}\left[\theta_{\alpha\alpha\alpha}\right]$ gives rise to a stiffness condition $\Delta t\leq C \left(\Delta s\right)^3$. The same analysis shows that the term $N\left(\alpha,t\right)$ is non-stiff. 

We identify that in Fourier space, the dominant term on the right hand side of the Eq.~(\ref{SSD}) diagonalizes and the equation becomes 
\begin{equation}
    \hat{\theta}_t = - \frac{\sigma \left|n\right|^3}{s_{\alpha}^3} \hat{\theta}\left(k,t\right)+\hat{N}\left(k,t\right). \label{SSDF}
\end{equation}
 We time-integrate the $\theta$-equation in Fourier space with a semi-implicit time-stepping algorithm \cite{HLS}. Using an integrating factor, we obtain
\begin{equation}
    \frac{d}{dt}\left(e^{-\frac{\sigma \left|n\right|^3}{s_{\alpha}^3}}\hat{\theta}_t\right) =e^{-\frac{\sigma \left|n\right|^3}{s_{\alpha}^3}}\hat{N}\left(k,t\right). \label{SSDF1}
\end{equation} 
Then, we use a second-order Adams-Bashforth (AB2) method to discretize Eq.~(\ref{SSDF1}) as
\begin{align}
    \hat{\theta}^{n+1}\left(k\right) &= e_k\left(t_n,t_{n+1}\right)\hat{\theta}^{n}\left(k\right)\nonumber\\
    & \quad +\frac{\Delta t}{2}\left(3e_k\left(t_n,t_{n+1}\right)\hat{N}^n\left(k\right)-e_k\left(t_{n-1},t_{n+1}\right)\hat{N}^{n-1}\left(k\right) \right),
\end{align}
where the subscript/superscript $n$ denotes numerical solution at $t=t_n$ and we define 
\begin{equation}
e_k\left(t_n,t_{n+1}\right) = \exp\left(-\sigma\left|k\right|^3\int_{t_n}^{t_{n+1}}\frac{dt}{s_{\alpha}^3\left(t\right)}\right).
\end{equation}
To evaluate the term $e_k\left(t_n,t_{n+1}\right)$, we first integrate the non-stiff Eq.~(\ref{lq}) using AB2 which gives
\begin{equation}
L^{n+1} = L^{n} +\frac{\Delta t}{2}\left(3M^n-M^{n-1} \right),
\end{equation}
with $M = -\frac{1}{2\pi}\int_{0}^{2\pi}V\left(\alpha,t\right)\theta_{\alpha} \, d\alpha$. Also, $s_{\alpha}=L/2\pi$, and we apply the trapezoidal rule to evaluate integrals in $e_k\left(t_n,t_{n+1}\right)$ and $e_k\left(t_{n-1},t_{n+1}\right)$ as
\begin{align}
\int_{t_n}^{t_{n+1}} \frac{dt}{s^3_{\alpha}\left(t\right)}&\approx \frac{\Delta t}{2}\left( \frac{1}{\left(s^n_{\alpha}\right)^3}+\frac{1}{\left(s^{n+1}_{\alpha}\right)^3}\right),\\
\int_{t_{n-1}}^{t_{n+1}} \frac{dt}{s^3_{\alpha}\left(t\right)}&\approx \Delta\left( 
\frac{1}{2\left(s^{n-1}_{\alpha}\right)^3}+
\frac{1}{\left(s^n_{\alpha}\right)^3}+\frac{1}{2\left(s^{n+1}_{\alpha}\right)^3}\right).
\end{align}
The AB2 method depends on two previous values, and therefore, we initiate the computation at time $t=0$ using Euler's method to obtain the relevant quantities at $t=\Delta t$. In the subsequent time-steps, we use the AB2 method as two previous time-step values are always known. The accumulation of noise is a problem \cite{JLL}; therefore, we employ a cutoff filter to prevent the accumulation of round-off error \cite{Krasny1} and a 25th-order Fourier filter to damp the higher, nonphysical modes and suppress the error due to aliasing.

\section{Numerical Results}\label{results}
In this section, we discuss the results of our numerical simulations. We first compare the results of nonlinear simulation with linear analysis and then demonstrate the spatio-temporal accuracy of our code. Finally, we compute several interesting cases where the domain $\Omega^{-}$ has different initial configuration. In all our simulations, we set the surface tension parameter to $\sigma=0.47$.

\subsection{Comparison of results of linear analysis and nonlinear simulation}
\begin{figure}[!ht]
    \centering
    \includegraphics[width=0.75\textwidth]{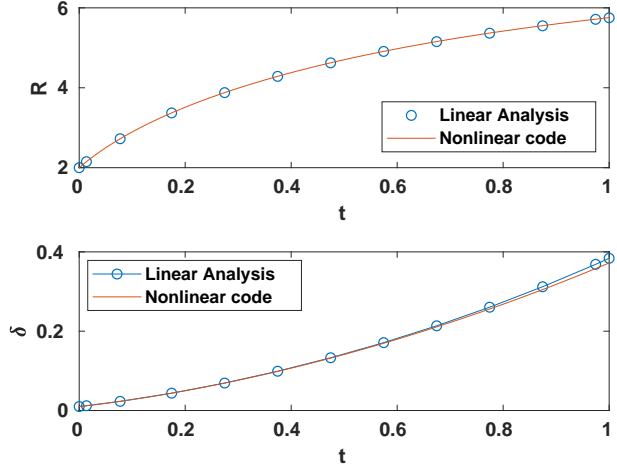}
    \caption{Comparison of results from the nonlinear simulation and the linear analysis for $R(t)$ and $\delta(t)$ against time. We choose $\sigma=0.47, R_{\infty}=10, N=1024$, and $\Delta t = 2\times 10^{-3}$ to obtain the match between the two setups and the simulation are stopped when the linear analysis results starts to over-predict the nonlinear results at $t_{\rm end}=1.0$.} 
    \label{fig:1c}
\end{figure}

The evolution of a perturbed circular interface is investigated, with the initial interface at $t=0$ given by
\begin{equation}
R+\delta \cos 4\theta =2 + 0.01\times \cos 4\theta \label{shape},
\end{equation}
and we choose $R_{\infty} = 10$. The simulation is carried out up to a time $t_{\rm end}=1.0$. Evolution of $R(t)$ and $\delta(t)$ against time are shown in Fig.~\ref{fig:1c}, using results from the nonlinear simulation and  the linear analysis (Eqs.~(\ref{R_t}) and(\ref{d_t})). The plots indicate excellent match between the two in the beginning thus validating our numerical methods. Once $\delta$ becomes large, we observe disagreement between the results of the linear analysis and the nonlinear simulation, especially in the evolution of $\delta$. It is evident from the plots that the linear system over-predicts the growth of the mode. This simulation confirms that the linear solution  holds for a short time span and the fully nonlinear simulation is needed to predict the evolution over a longer time.

\begin{figure}[!ht]
    \centering
    \includegraphics[width=0.75\textwidth]{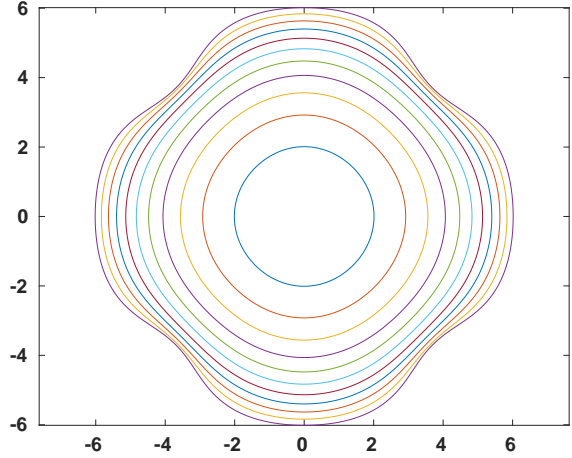}
    \caption{Time evolution of the interface}
    \label{fig:1d}
\end{figure}

Fig.~\ref{fig:1d} shows the evolution of the interface, where the innermost contour corresponds to the shape at $t=0$. For all simulations up to this point, we used a GMRES tolerance of $\epsilon = 10^{-10}$. The filters are also set to this tolerance.

\subsection{Spatio-Temporal convergence}

\begin{figure}[ht!]
\centering
\begin{subfigure}[Spectral accuracy]{
\centering
\includegraphics[width=0.8\textwidth]{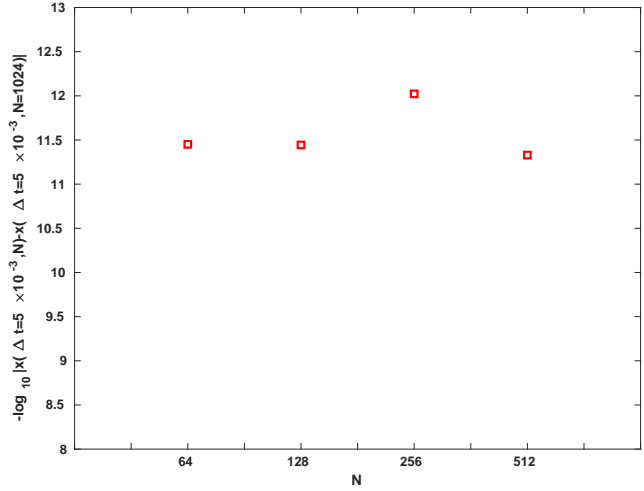}
\label{fig:1a}}
\end{subfigure}%
\begin{subfigure}[Temporal accuracy]{
\centering
\includegraphics[width=0.8\textwidth]{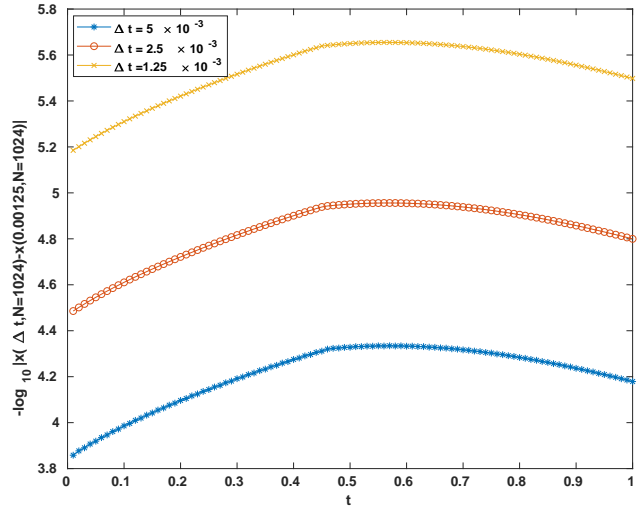}
\label{fig:1b}}
\end{subfigure}
\caption{Demonstration of spectral accuracy and second-order convergence in time of the nonlinear simulation.} 
\end{figure}

Figs.~\ref{fig:1a} and \ref{fig:1b} show the spatio-temporal accuracy of our numerical simulation using initial shape defined in Eq.~(\ref{shape}) and with other parameters unchanged. Note that our numerical method is spectrally accurate is space and second-order accurate in time.
In Fig.~\ref{fig:1a}, we demonstrate the spectral accuracy of our code by plotting the maximum of $$-\log_{10}\left|\text{x}\left(t,N\right)-\text{x}\left(t,N_f=1024\right)\right|$$
for values $N=64, 128, 256$, and $512$ at time $t_{\rm end}=1.0$. $\Delta t = 5\times 10^{-3}$ is chosen so that the results are very accurate in time. Observe that even with $N=64$, the results match up to $10^{-11}$. This indicates very a rapid decay of error with $N$ and confirms the spectral accuracy of our code.

In Fig.~\ref{fig:1b}, we plot the maximum of $-\log_{10}|\text{x}\left(\Delta t,N\right)-\text{x}\left(5 \times 10^{-4}, N\right)|$ for $N=1024$ and three values of $\Delta t = 5\times 10^{-3}$, $2.5\times 10^{-3}$, and $\Delta t = 1.25 \times 10^{-3}$ until the time $t_{\rm end}=1$. The distance between the lines is 0.6, indicating second-order convergence. 
We deliberately choose large $N$ during temporal convergence study 
to ensure high accuracy in space such that the space discretization error does not interfere with the error due to time discretization.

\subsection{Simulation of different steady state configurations}
In this section, we show different steady state configurations starting with various initial conditions. We set the GMRES tolerance to $\epsilon= 10^{-8}, N=512$, filter tolerance to $10^{-10}$, and  $\Delta t = 5\times10^{-4}$ unless stated otherwise. We found that the relaxed tolerance does not interfere with the accuracy of simulation, but a stricter temporal resolution helps improve convergence. We further found that $N=512$ is enough for space resolution throughout the simulation as the morphologies are not complicated. All simulations except the last one are performed using an Intel(R) Core(TM) i5-7200U processor with maximum clock speed @ 2.50GHz and in a laptop with 8 GBs of RAM space. The last simulation with 12 regions was carried out on a desktop machine with Intel(R) Core(TM) i9-10900 processor with maximum clock speed @ 2.80GHz and 64 GB RAM. 

In all our simulations, we maintain the following protocol: We start the simulation under transient conditions where the system is driven by the flux given in Eq.~(\ref{bc3new}). Once the right hand side of the equation is less than a tolerance value of 0.001, we set the flux forcefully to zero. We do this because the flux goes to zero only as $t\rightarrow \infty$ but, for all practical purposes,  can be neglected when it goes below the small tolerance we set. Once that happens, 
the system  moves into the zero-flux regime or the relaxation phase and we observe the evolution for sufficiently long time to investigate the domain configurations in the steady-state. We stop the simulation at $t_{\rm end}=25$ if it does not stop earlier due to a topological singularity showing up in the system. 
In time plots, we always use semilog in the $x$-axis. 

\begin{figure}[ht!]
\centering
\begin{subfigure}[$t-$Maximum normal velocity plot]{
\centering
\includegraphics[width=0.95\textwidth,trim = {0 2.0cm 0 4.0cm}]{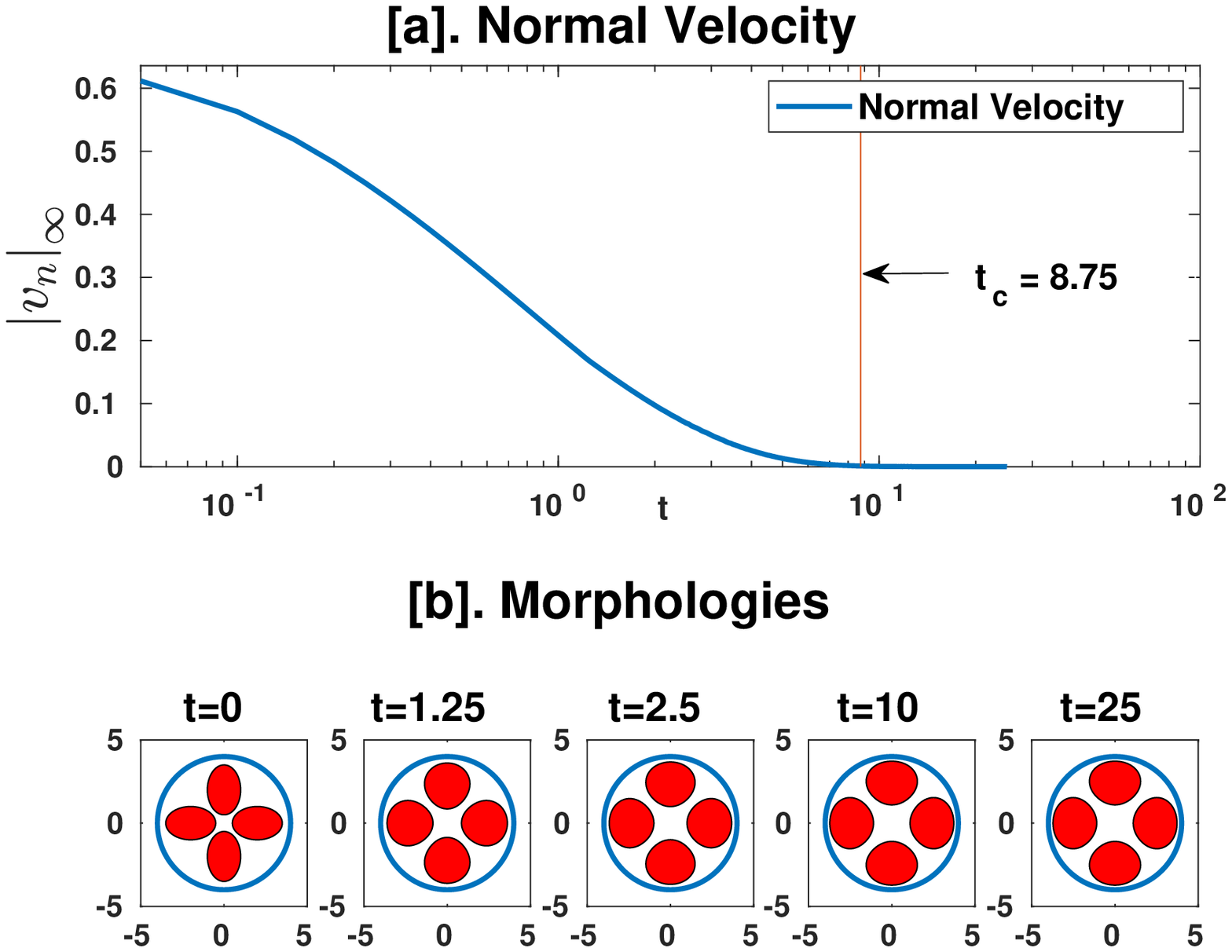}
\label{four_fig1A} }
\end{subfigure}
\begin{subfigure}[$t-w_{\infty}$ plot]{
\centering
\includegraphics[width=0.46\textwidth]{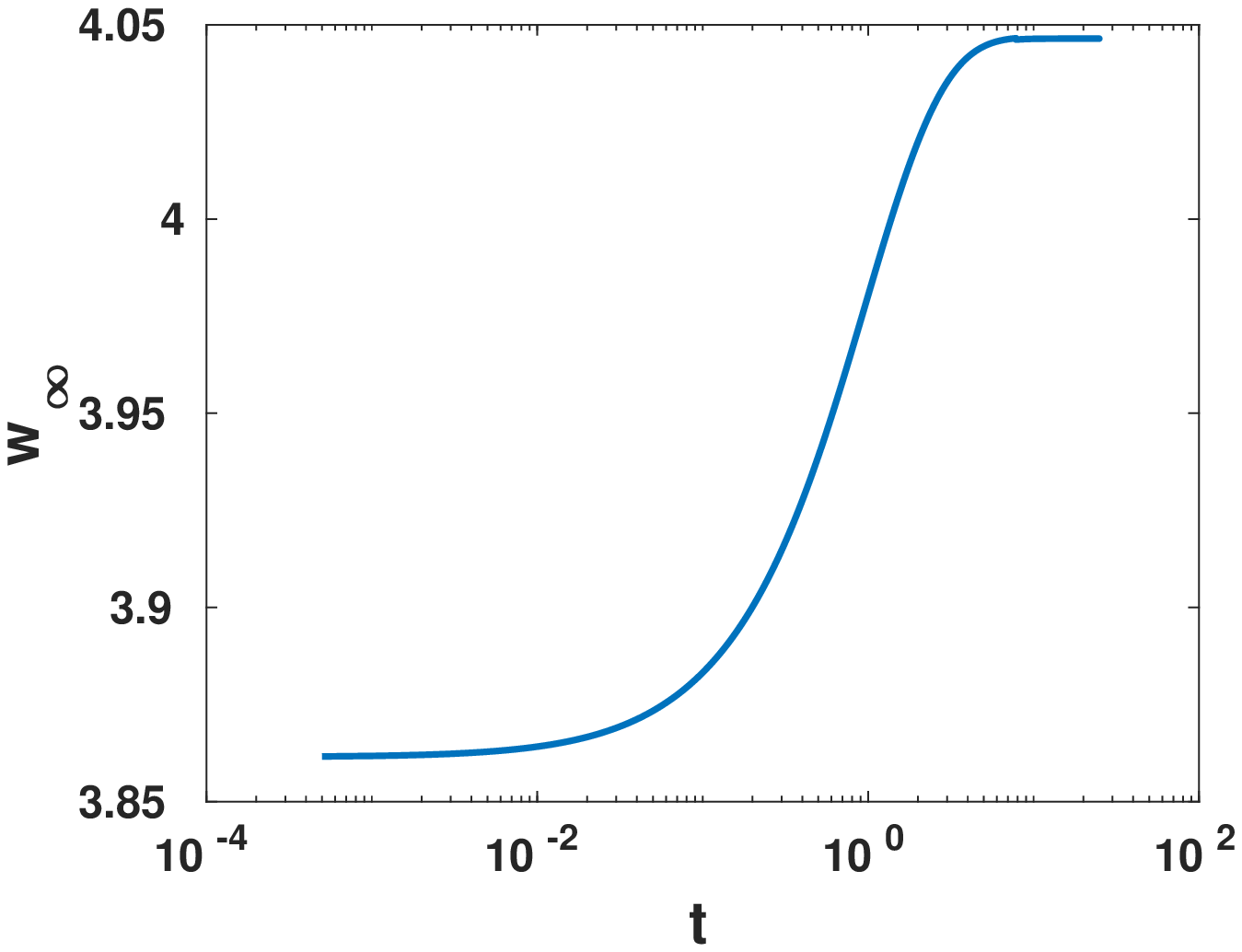}
\label{four_fig1B} }
\end{subfigure}
\begin{subfigure}[$t-s_{\alpha}$ plot]{
\centering
\includegraphics[width=0.46\textwidth]{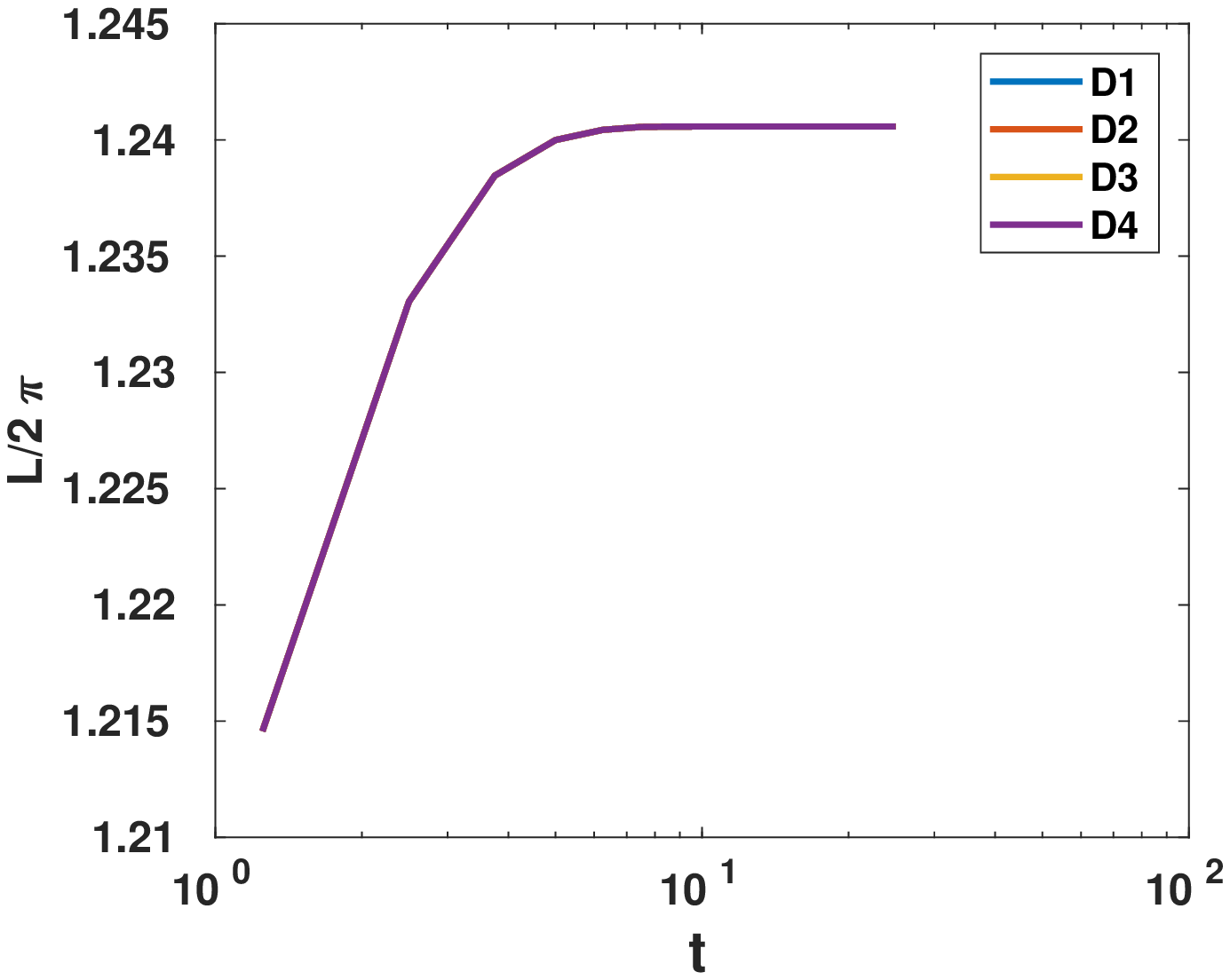}
\label{four_fig1C} }
\end{subfigure}
\caption{Time evolution of 4 elliptic regions with semi-axes $a=1.5$ and $b=1.0$. The other parameters are $R_{\infty}=4$ and surface tension $\sigma=0.47$. The system enters equilibrium at $t_{eq}=8.75.$ Centroids of the domains D1, D2, D3, and D4 are at 
$\left(2,0\right)$, 
$\left(0,2\right)$, 
$\left(-2,0\right)$, and
$\left(0,-2\right)$ at $t=0$, respectively.}
\end{figure}

First, we perform a simulation using a four-domain configuration and display the results of various important parameters of the simulation in Fig.~\ref{four_fig1A}, Fig.~\ref{four_fig1B}, and Fig.~\ref{four_fig1C}. The domains at $t=0$ are elliptic in shape and we have one domain each along the positive and negative $x$- and $y$-axes. The major and minor axes of the ellipses are set to the values $a=1.5$ and $b=1.0$. We set $R_{\infty}=4$ and the centroids of the domains are at $
\left(2,0\right), \left(0,2\right), 
\left(-2,0\right),
\left(0,-2\right).$
We denote these domains by D1, D2, D3, and D4, respectively. The initial configuration (lower left panel of Fig.~\ref{four_fig1A}) is symmetric about the $x$- and $y$-axes. It also has certain rotational symmetries. The governing equations demand that these symmetries should be preserved at all later times and we find that this is indeed true for our simulation. 

With this configuration, we find that the changes are rapid at the beginning. The outer parts of the ellipses bulge out and align themselves along the boundary perhaps because more space is available towards the outer region as compared to region near the center, and by time $t=2.5$, the  shapes no more resemble ellipses.  
The system enters equilibrium configuration at $t_c=8.75$ when the flux approaches zero. To understand more about this phase, we refer to the plot of the maximum interfacial velocity $\max \left\|v\right\| = \left\|v\right\|_{\infty}$ (top panel of Fig.~\ref{four_fig1A}) where the maximum is taken over all marker points over all interfaces. It is observed in this plot that the velocity decreases monotonically to zero, and close to $t_c$, the maximum magnitude of the velocity $\max \left\|v\right\| = \left\|v\right\|_{\infty}$ is negligible. Therefore, the system configuration changes very little in the relaxation phase. This is confirmed by comparing the plots of the configuration (lower panels of Fig.~\ref{four_fig1A}), in which the changes after $t=2.5$ are small. At $t=t_{\rm end}$, we find that the domains lose their elliptic form and are approximately circular.

The evolution of two additional quantities, the arc length parameter $s_{\alpha}=L/2\pi$ for each interface, and the far-field function value $w_{\infty}$, are shown in Fig.~\ref{four_fig1C} and  Fig.~\ref{four_fig1B}, respectively. Because of the symmetry, all four curves are on the top of each other in Fig.~\ref{four_fig1C}. The far-field flux is flat at the beginning but eventually changes rapidly before entering the relaxation phase, giving it the shape of a sigmoid curve. The change in arc length parameter is rapid at the beginning but this curve flattens very quickly once the system enters the relaxation phase.

\begin{figure}[ht!]
\centering
\begin{subfigure}[$t-$Maximum normal velocity plot]{
\centering
\includegraphics[width=0.95\textwidth,trim = {0 1.2cm 0 4.0cm}]{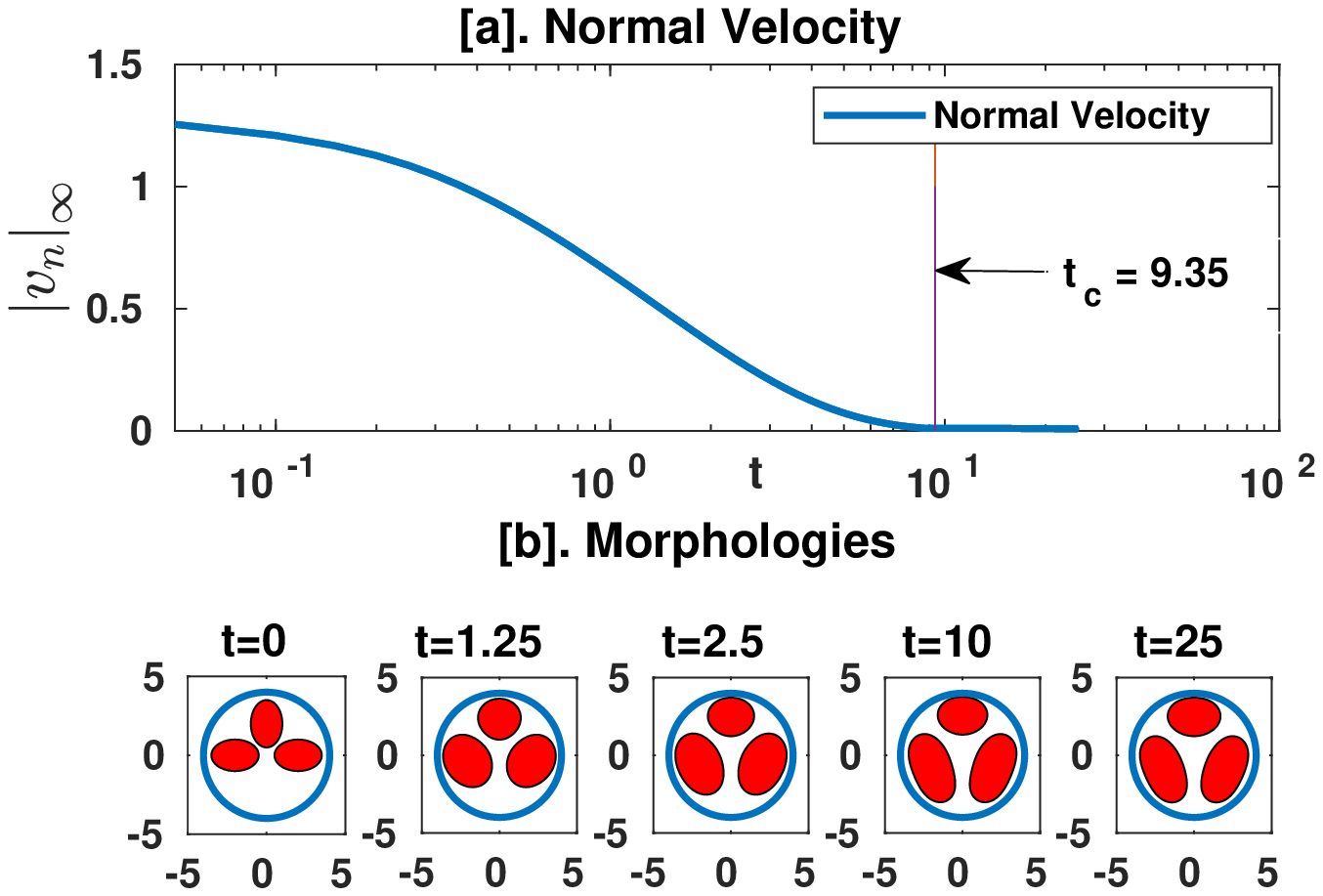}
\label{three_fig1A} }
\end{subfigure}
\begin{subfigure}[$t-w_{\infty}$ plot]{
\centering
\includegraphics[width=0.46\textwidth]{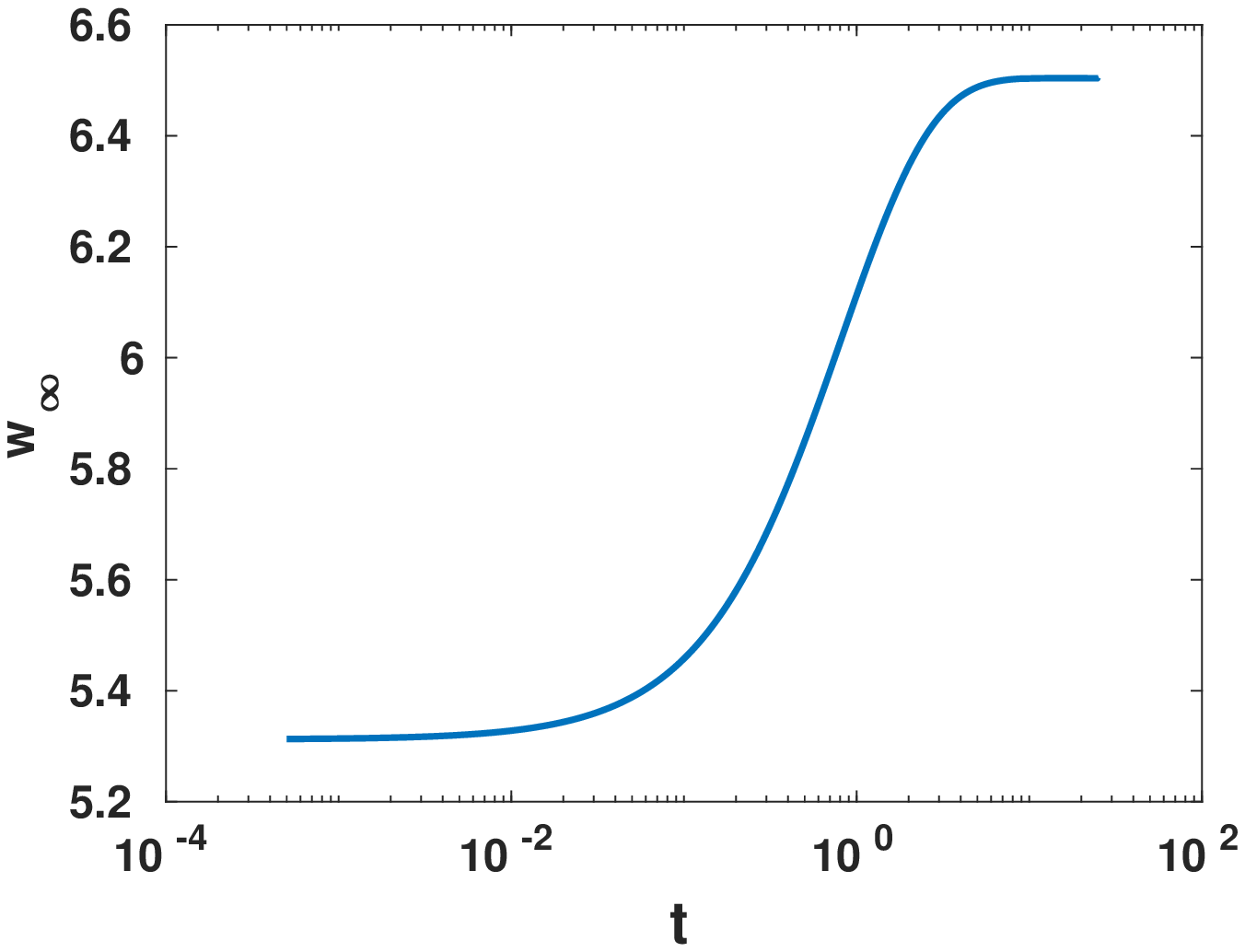}
\label{three_fig1B} }
\end{subfigure}
\begin{subfigure}[$t-s_{\alpha}$ plot]{
\centering
\includegraphics[width=0.46\textwidth]{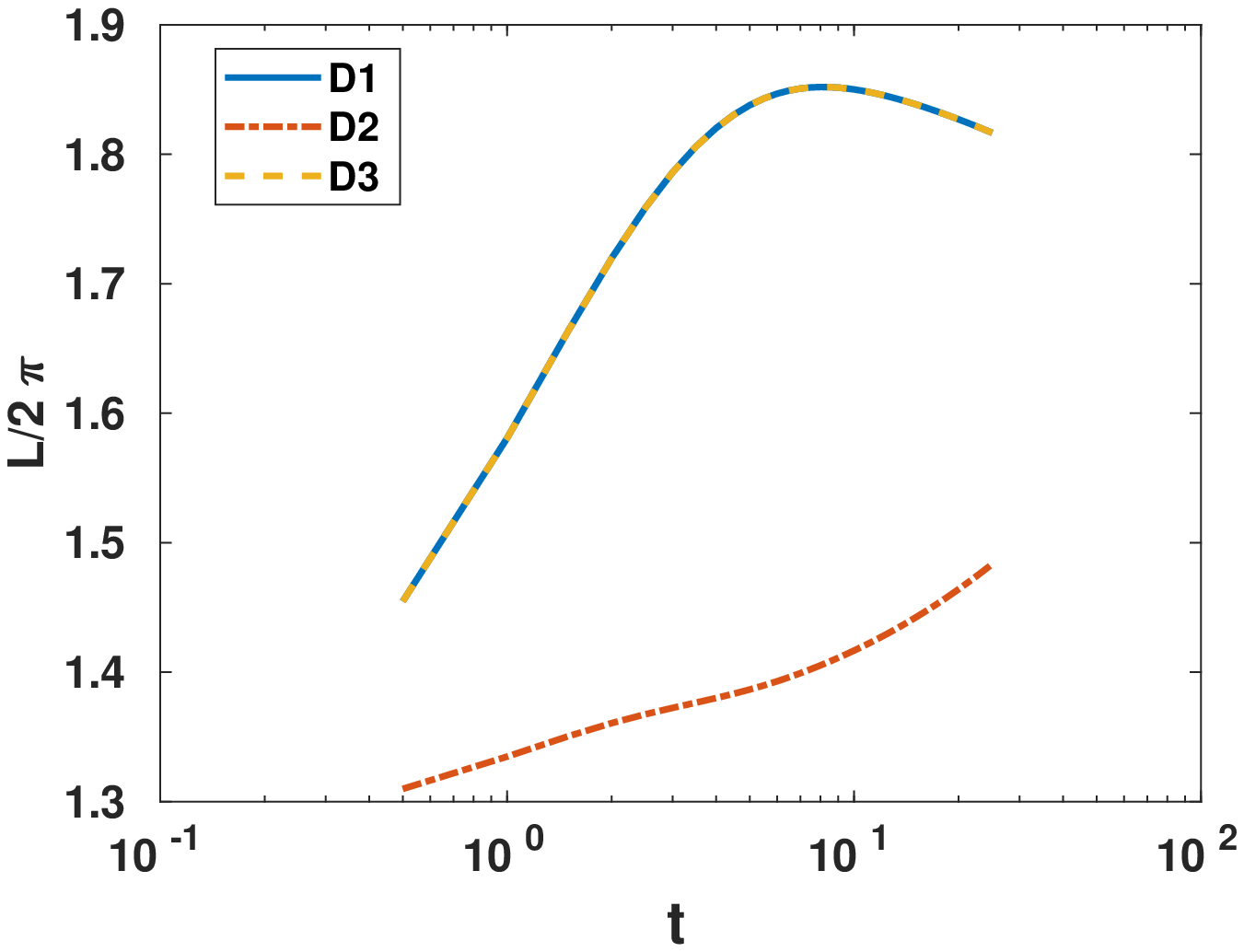}
\label{three_fig1C} }
\end{subfigure}
\caption{Time evolution of 3 elliptic regions with semi-axes $a=1.5$ and $b=1$. We set $R_{\infty}=4$. The system enters the equilibrium phase at $t_{c}=9.35.$ Centroids of the domains D1, D2, and D3 are at $
\left(2,0\right)$, 
$\left(0,2\right)$, and
$\left(-2,0\right)$ at $t=0$, respectively.}
\end{figure}

Next, we consider a simulation with three domains. We do this by removing one particle from the previous configuration.
In Fig.~\ref{three_fig1A}, the initial configuration  is symmetric about both the $x$- and $y$-axis. We start with elliptic particles with semi-axes dimensions of $a=1.5$ and $b=1.0$, and with their centroids at $\left(2,0\right), \left(0,2\right)$ and $\left(-2,0\right)$. We label these regions D1, D2, and D3, respectively.
The radius of the far-field boundary is $R_{\infty}=4$.

We observe that the domains D1 and D3, originally aligned along positive and negative $x$-direction respectively, rotate quickly, by almost 45 degrees. By $t=1.25$, significant rotation occurs and it continues further, even as the angular speed slows down. 
The domain D2 shrinks in the $y$-direction and grows in the $x$-direction. After sufficient time, this domain is ellipse-like with major axis in $x$-direction and minor axis along $y$-axis.

An interesting point is the difference in the area occupied by each domain as the time progresses. The area of the domains are equal in the beginning. As the simulation progresses, all regions grow in size, with region D2 growing slower the other two particles. This is  prominent during the early stages of evolution. However, the area of D2 increases somewhat faster during the later stages of evolution (after $t=10$), and eventually, the ratio of the arc length parameters of D2 and D1/D3 is  approximately $1.2$. The flux approaches zero at approximately $t_c=9.35$. 

\begin{figure}[ht!]
\centering
\begin{subfigure}[$t-$Maximum normal velocity plot]{
\centering
\includegraphics[width=0.95\textwidth,trim = {0 1.2cm 0 4.0cm}]{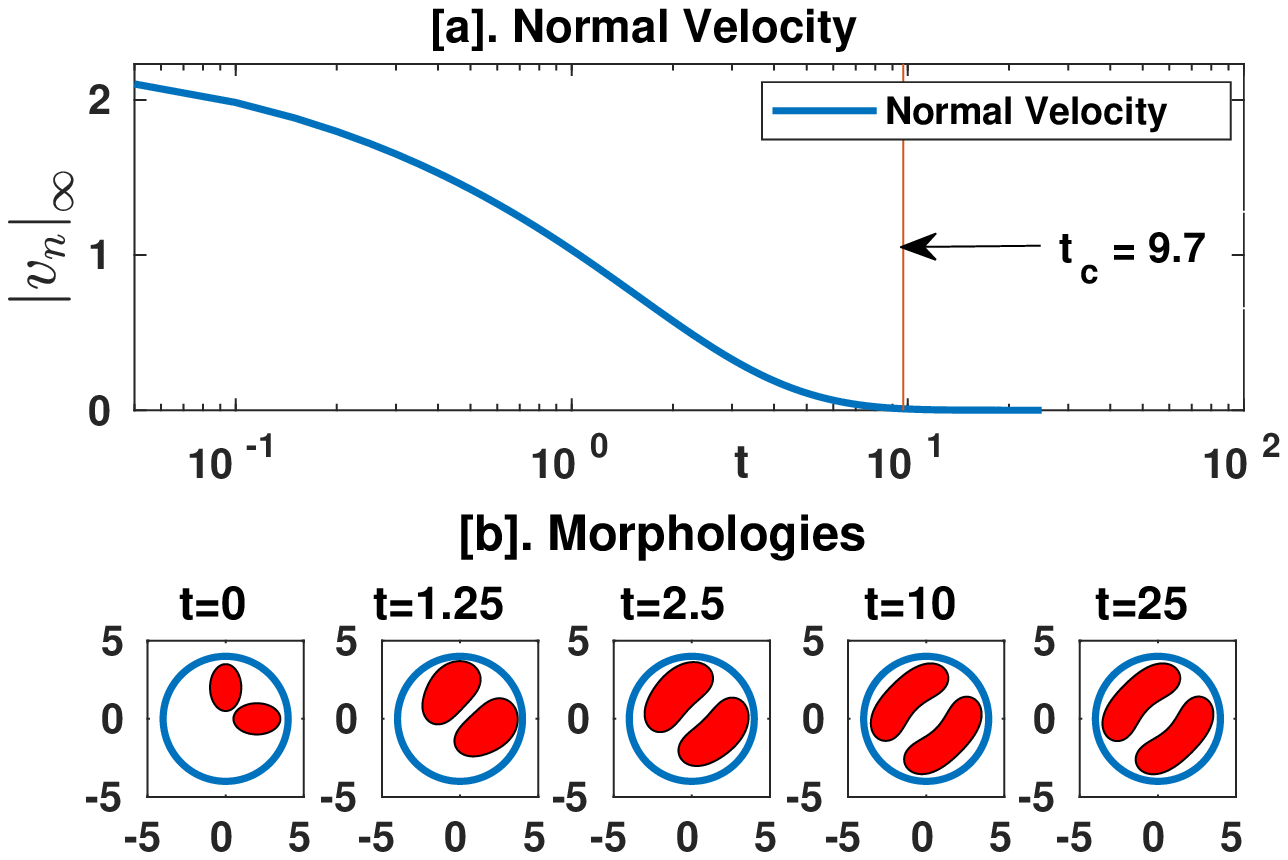}
\label{two_fig1} }
\end{subfigure}
\begin{subfigure}[$t-w_{\infty}$ plot]{
\centering
\includegraphics[width=0.46\textwidth]{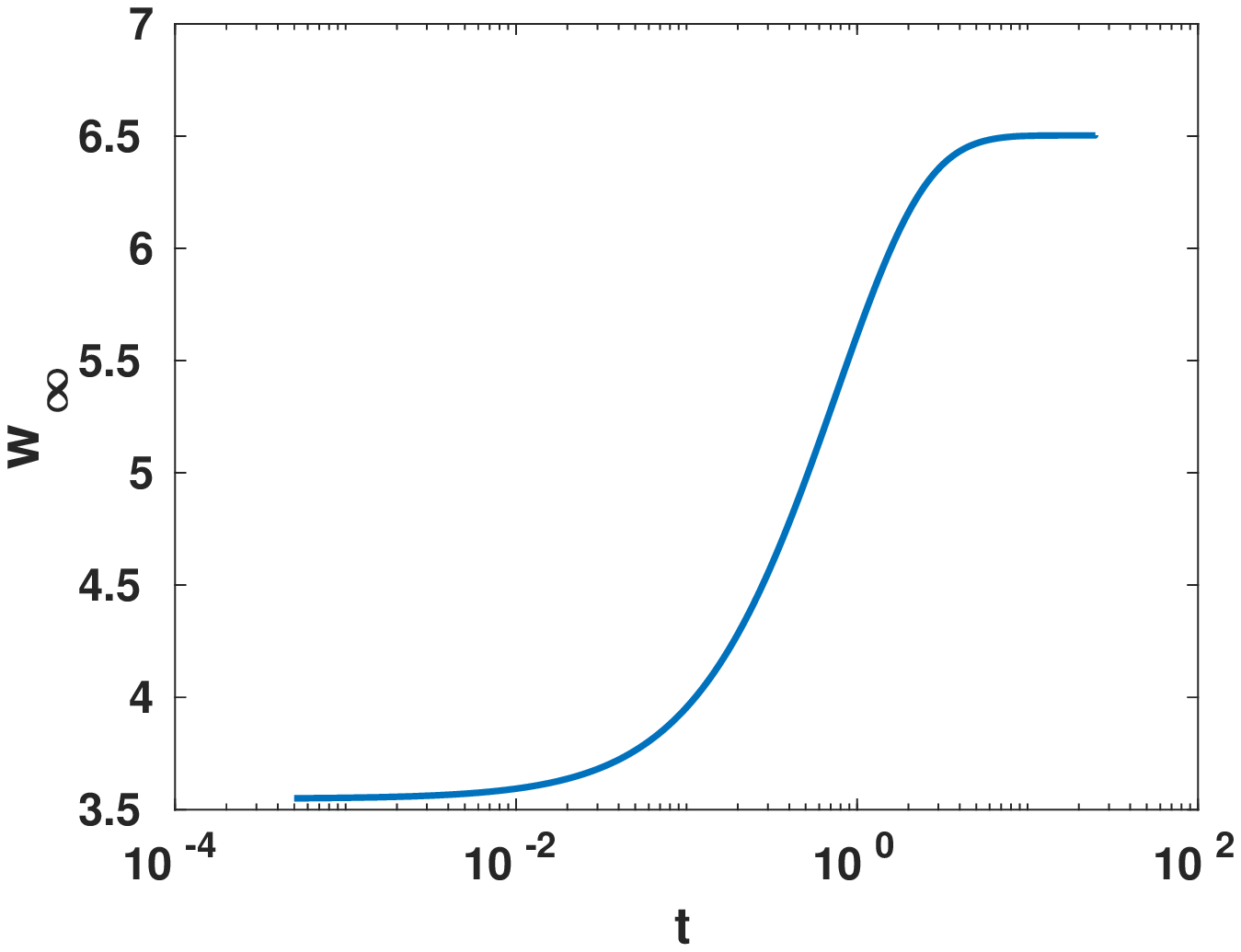}
\label{two_fig2} }
\end{subfigure}
\begin{subfigure}[$t-s_{\alpha}$ plot]{
\centering
\includegraphics[width=0.46\textwidth]{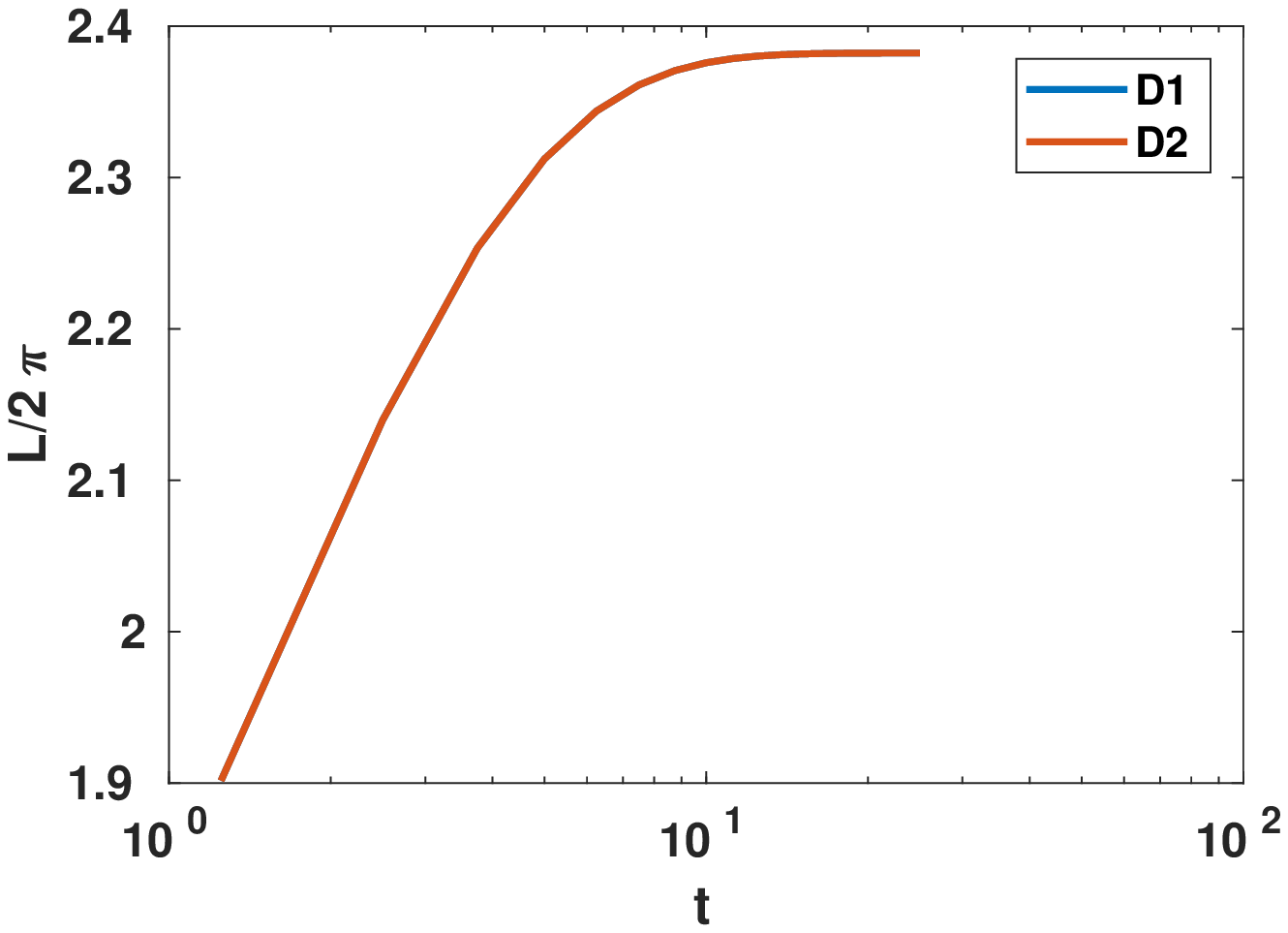}
\label{two_fig3} }
\end{subfigure}
\caption{Time evolution of 2 elliptic regions with semi-axes $a=1.5$ and $b=1$. We set $R_{\infty}=4$. The system enters equilibrium at $t_{c}=9.7.$ Centroids of the domains D1 and D2 are at $\left(2,0\right)$ and $\left(0,2\right)$ at $t=0$, respectively.}
\end{figure}

Figs.~\ref{two_fig1},~\ref{two_fig2}, and~\ref{two_fig3} show simulation results corresponding to two elliptic phase domains. The domains are aligned along the $x$- and $y$-axes with semi-axes dimensions $a=1.5$ and $b=1.0$. We set $R_{\infty}=4$. 
The centroid of the phase domain with major axes along $x$-direction is at $
\left(2,0\right)$, and the other one is located at $\left(0,2\right)$. This configuration is symmetric about the line $y=x$. The domains undergo rotation during evolution, aligning themselves along the line $x=y$ and growing in size during the alignment process due to a positive flux. The particle shapes are convex towards the boundary $\partial \Omega$ while they are concave in the inner region. As with the simulation with four and three domains, the graph of $w_{\infty}$ has a sigmoid shape. 

\subsection{Domain shrinkage}
\begin{figure}[ht!]
\centering
\begin{subfigure}[$t-$Maximum normal velocity plot]{
\centering
\includegraphics[width=0.95\textwidth,trim = {0 2.0cm 0 4.0cm}]{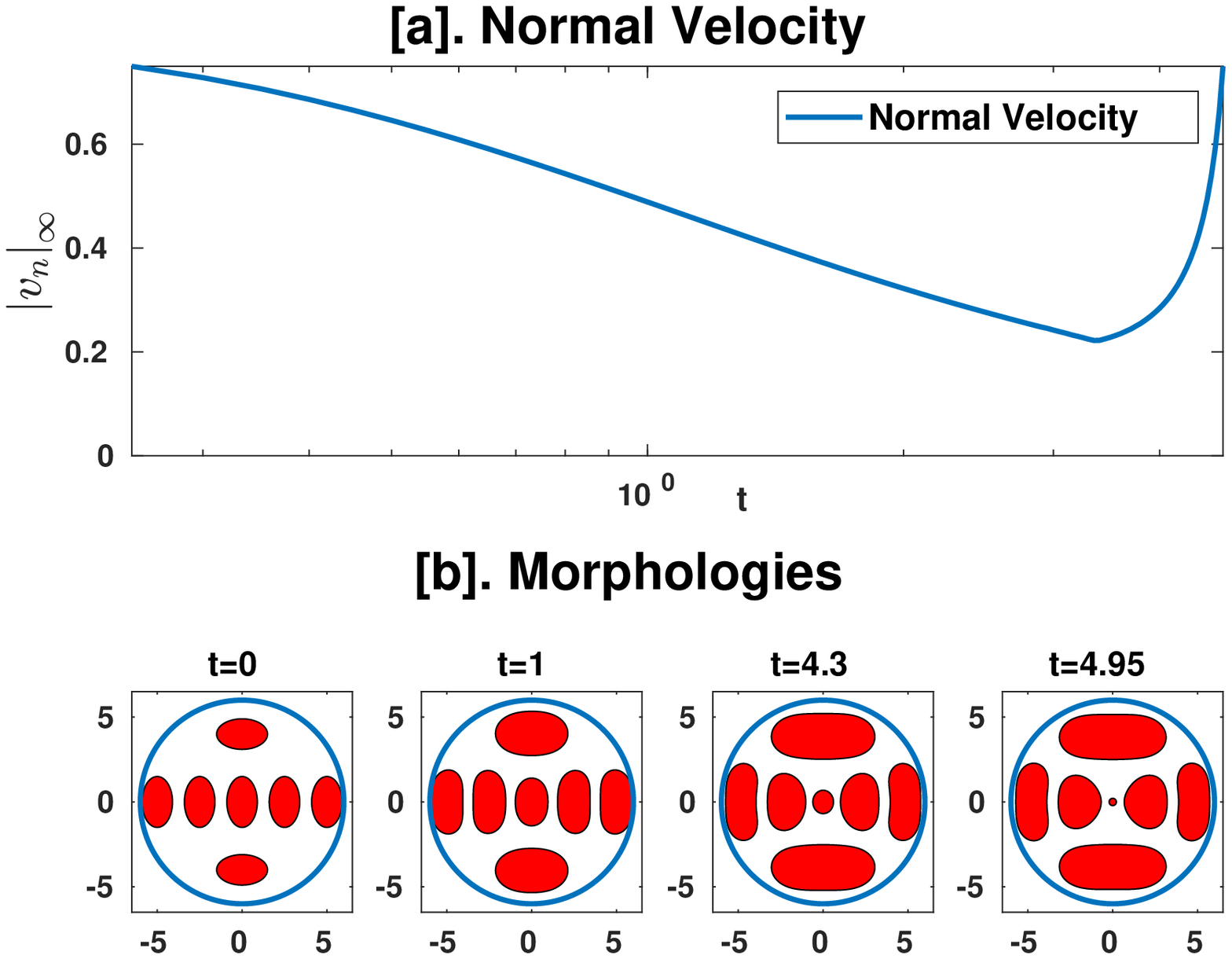}
\label{seven_fig1} }
\end{subfigure}
\begin{subfigure}[$t-w_{\infty}$ plot]{
\centering
\includegraphics[width=0.46\textwidth]{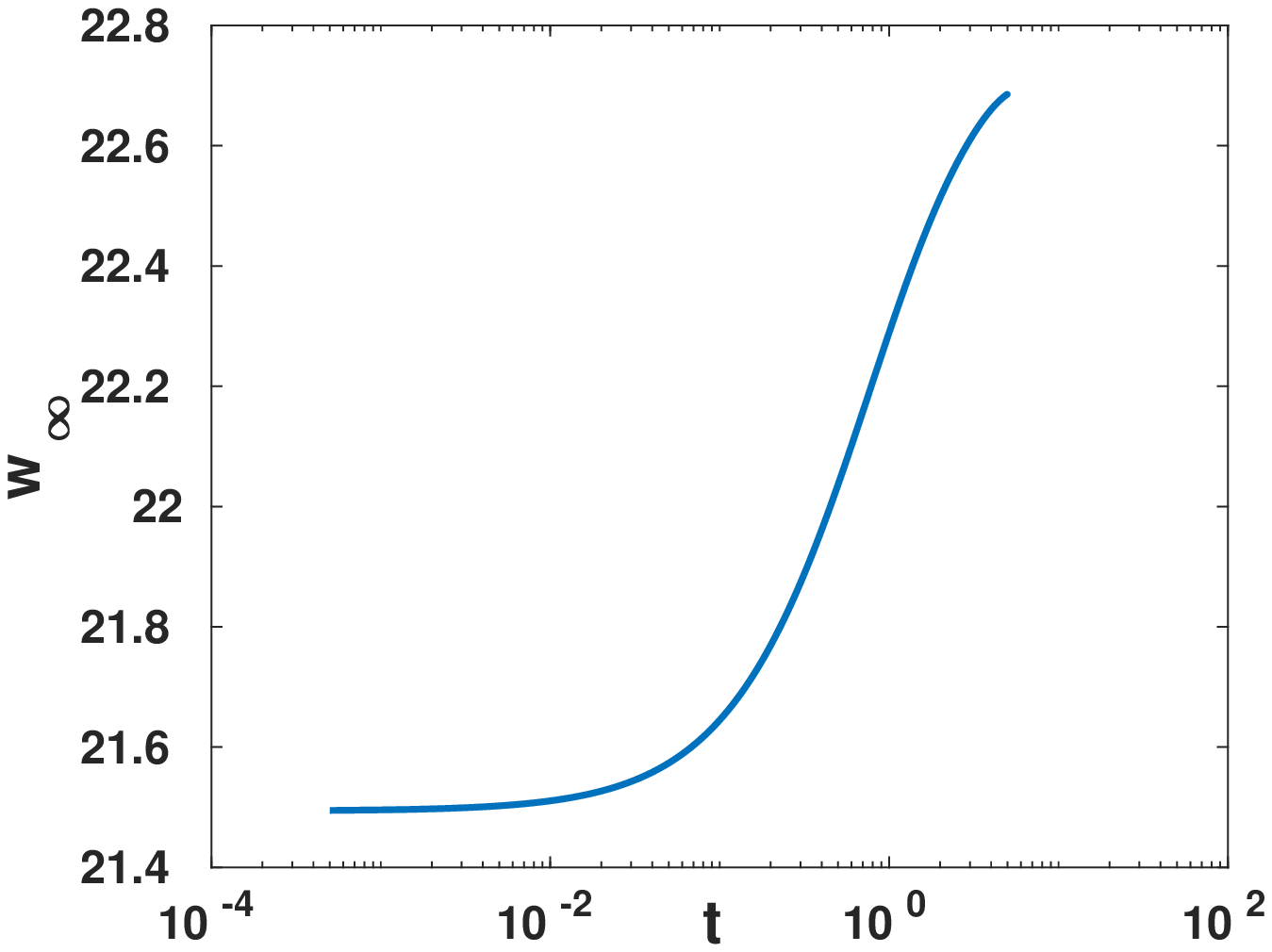}
\label{seven_fig2} }
\end{subfigure}
\begin{subfigure}[$t-s_{\alpha}$ plot]{
\centering
\includegraphics[width=0.46\textwidth]{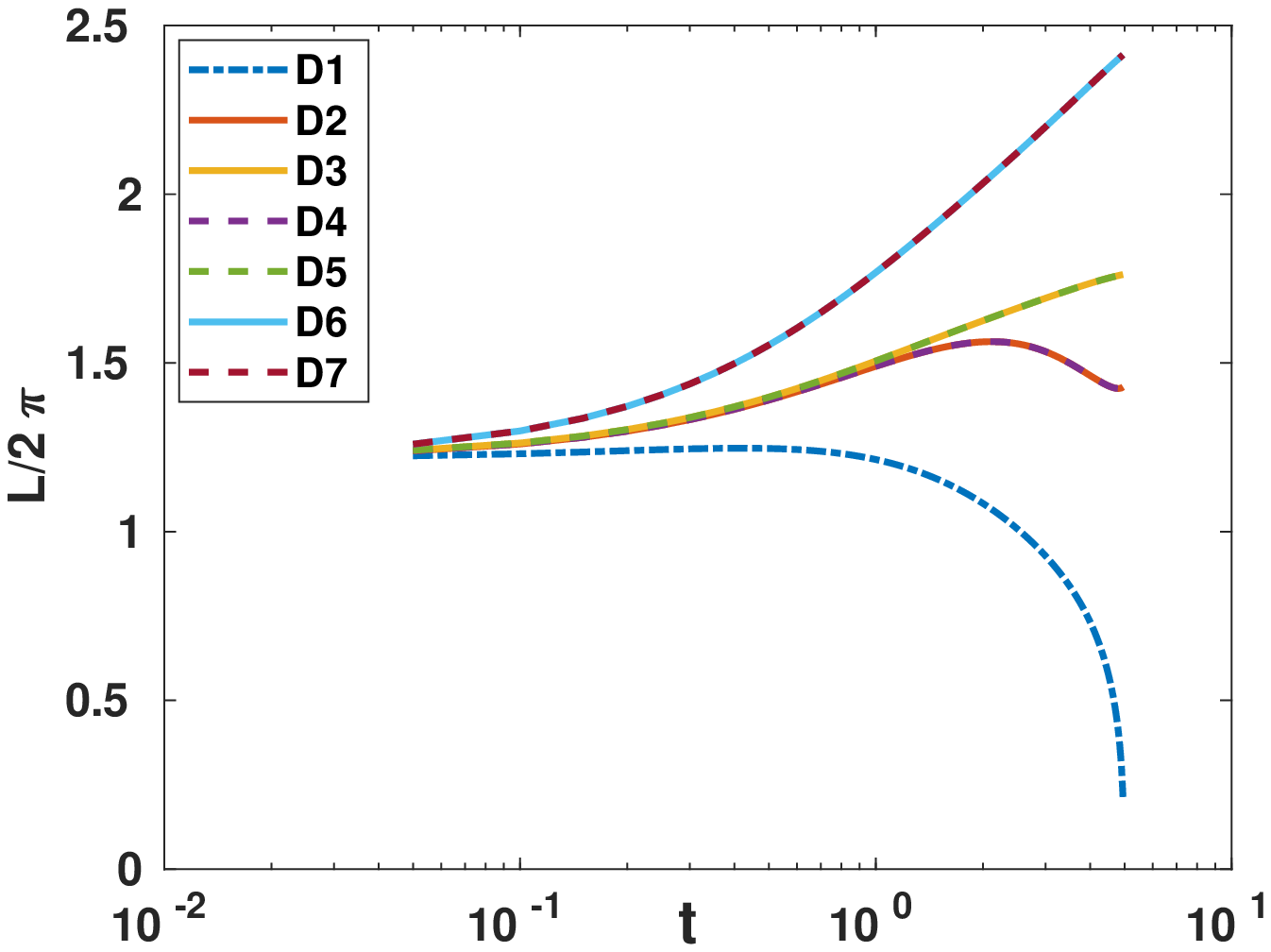}
\label{seven_fig3} }
\end{subfigure}
\caption{Time evolution of 7 elliptic regions with semi-axes $a=1.5$ and $b=0.9$. We set $R_{\infty}=6$. 
Centroids of the domains D1 to D7 are at 
$\left(0,0\right)$, 
$\left(2.5,0\right)$, 
$\left(5,0\right)$, 
$\left(-2.5,0\right)$, 
$\left(-5,0\right)$, 
$\left(0,4\right)$, 
and $\left(0,-4\right)$ at $t=0$, respectively.}
\end{figure}

Figs.~\ref{seven_fig1}, \ref{seven_fig2}, and \ref{seven_fig3} are results of simulations with seven elliptic domains. The centroids of the domains are at $
\left(0,0\right)$, $\left(2.5,0\right)$, 
$\left(5,0\right)$, $\left(-2.5,0\right)$, $\left(-5,0\right)$, $\left(0,4\right)$, and $\left(0,-4\right)$ with major axis $a=1.5$ and minor axis $b=0.9$. We denote these domains by D1 to D7, respectively. The outer boundary is at $R_{\infty}=6$. The configuration is symmetric about the $x$- and $y$-axes and has a rotational symmetry of 180 degrees. 

The evolution of this seven-domain configuration reveals a number of interesting aspects. Most notable of these is the shrinkage and gradual disappearance of the domain D1. All domains at $t=0$ have the same area but as time progresses, D1 shrinks. In the beginning, the area shrinks slowly but later the shrinking process speeds up. We note that near the singularity, around $t=4.75$, the code crashes and the results may not be very accurate. This is evident in the velocity plot where the maximum normal velocity decays at first and then increases very rapidly towards the end. Thus, our fixed time-steps may not capture the results towards the end of the simulation very well. 
The domains D6 and D7 are the ones that grow the most in the process. After these, the next largest growths are seen for D3 and D5, and then for D2 and D4. The arc lengths of the domains D2 and D4 display non-monotonic behavior with time. 

As a related phenomenon, we mention here the problem of particle coarsening~\cite{TAV1, TAV2, JLL} in alloy formation where, once the system enter the relaxation phase, the  phase-domains may undergo  topological changes. The domains tend to acquire compact shapes owing to the minimum surface energy requirements, and in the process, large domains try to grow at the expense of smaller regions. In this simulation, we find  results analogous to that.

\begin{figure}[ht!]
\centering
\begin{subfigure}[$t-$Maximum normal velocity plot]{
\centering
\includegraphics[width=0.95\textwidth,trim = {0 2.0cm 0 4.0cm}]{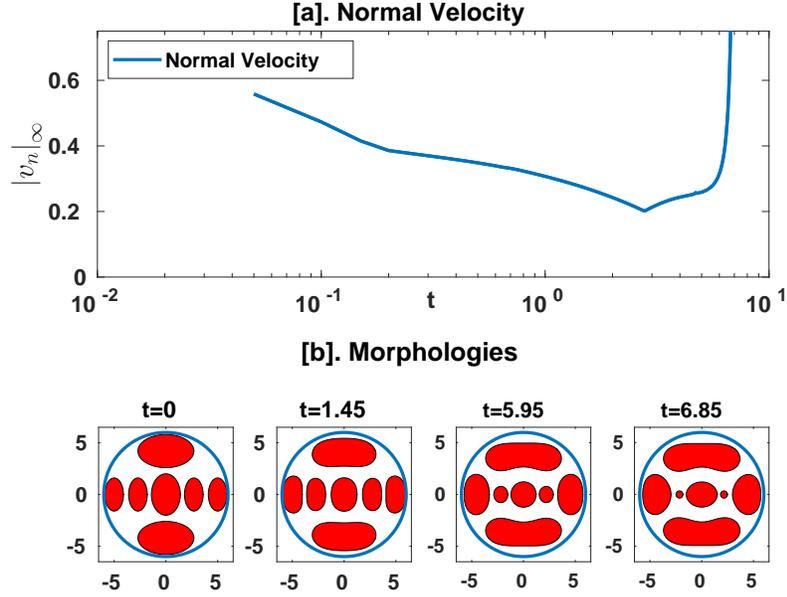}
\label{seven_fig2B}}
\end{subfigure}
\begin{subfigure}[$t-w_{\infty}$ plot]{
\centering
\includegraphics[width=0.46\textwidth]{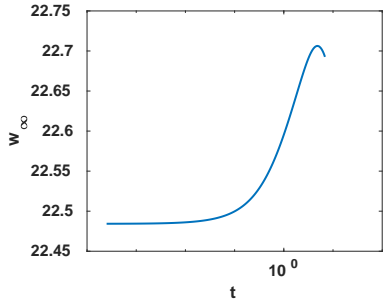}
\label{seven_fig2C} }
\end{subfigure}
\begin{subfigure}[$t-s_{\alpha}$ plot]{
\centering
\includegraphics[width=0.46\textwidth]{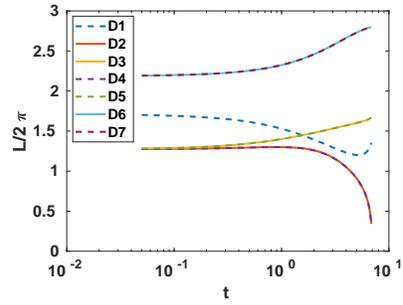}
\label{seven_fig2D} }
\end{subfigure}
\caption{Time evolution of 7 elliptic regions D1 to D7 with centroids at $\left(0,0\right)$, 
$\left(2.7,0\right)$, 
$\left(5,0\right)$, 
$\left(-2.7,0\right)$, 
$\left(-5,0\right)$, 
$\left(0,4\right)$, 
and $\left(0,-4\right)$ at $t=0$, respectively. The domain D1 has major and minor axes $a=2.0$ and $b=1.4$, domains D2 to D5 have major and minor axes $a=1.6$ and $b=0.9$, and domains D6 and D7 have major and minor axes $a=2.7$ and $b=1.6$. We set $R_{\infty}=6$.}
\end{figure}

Figs.~\ref{seven_fig2B}, \ref{seven_fig2C}, and \ref{seven_fig2D}, show results of a different seven-domain configuration. In this simulation, the regions D1 to D7 have their centroids at $\left(0,0\right)$, 
$\left(2.7,0\right)$, 
$\left(5,0\right)$, 
$\left(-2.7,0\right)$, 
$\left(-5,0\right)$, 
$\left(0,4.2\right)$, 
and $\left(0,-4.2\right)$ at $t=0$, respectively. The domain D1 has major and minor axes $a=2.0$ and $b=1.4$, domains D2 to D5 have major and minor axes $a=1.6$ and $b=0.9$, and domains D6 and D7 have major and minor axes $a=2.7$ and $b=1.6$. The areas of domains D1, D2, and D4 all decrease with the domains D2 and D4 shrinking faster than D1. This is in contrast with our previous simulation where D1 decreases fastest. Eventually D1 survives, but D2 and D4 disappears. Also the orientation of D1 changes, at time $t=0$ the major axis of D1 is aligned in $y$-direction, in an intermediate stage it is circular but towards the end it regains its elliptic shape to a certain extent and the major axis is in $x$-direction.  In this simulation we use a time step $\Delta t = 2.5\times 10^{-4}$, unlike in other simulations, as the reduced time step improves convergence.

\subsection{Simulation with large number of domains}
\begin{figure}[ht!]
\centering
\begin{subfigure}[$t-$Maximum normal velocity plot]{
\centering
\includegraphics[width=0.95\textwidth,trim = {0 2.0cm 0 4.0cm}]{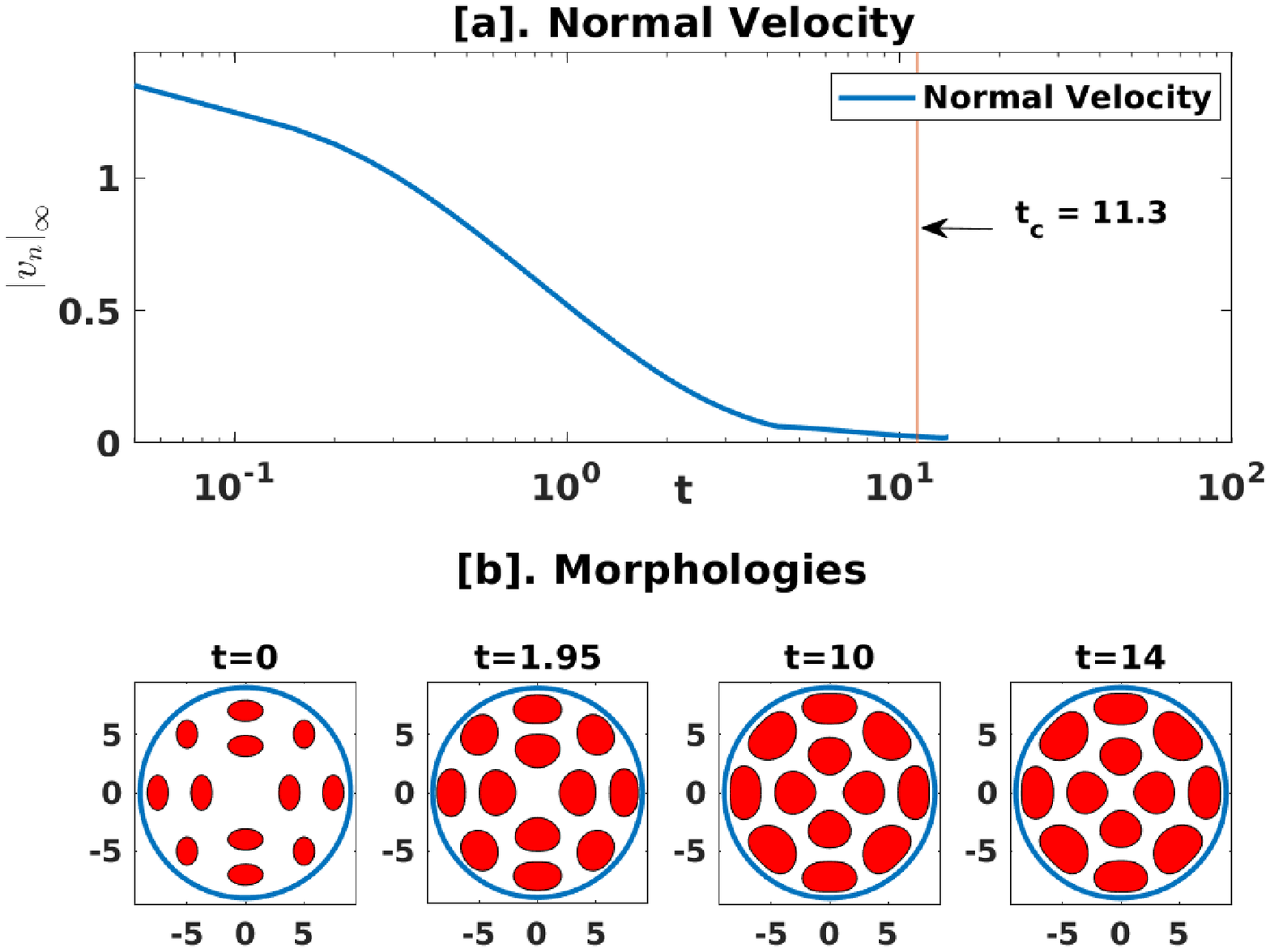}
\label{twelve_fig2B} }
\end{subfigure}
\begin{subfigure}[$t-w_{\infty}$ plot]{
\centering
\includegraphics[width=0.46\textwidth]{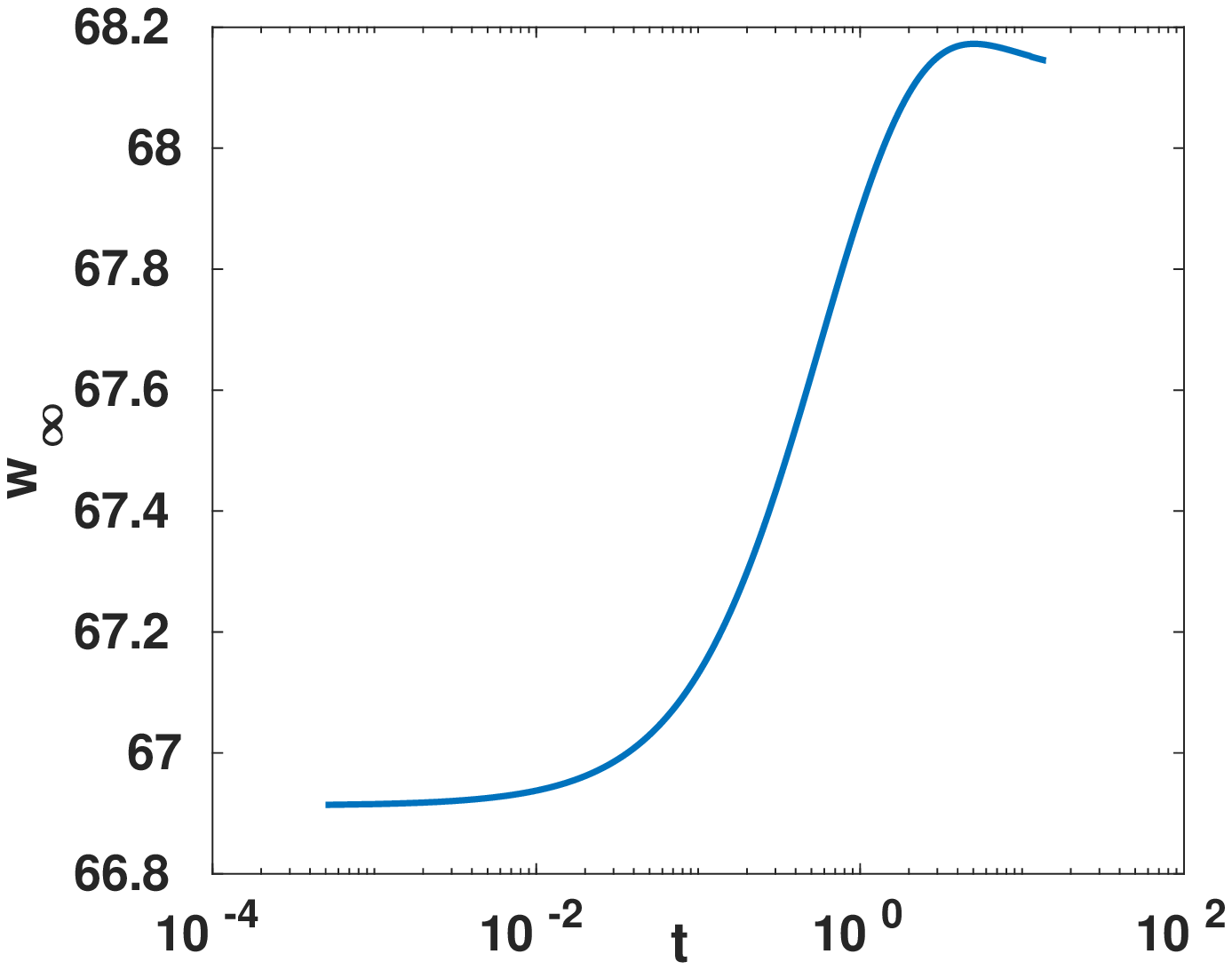}
\label{twelve_fig2C} }
\end{subfigure}
\begin{subfigure}[$t-s_{\alpha}$ plot]{
\centering
\includegraphics[width=0.46\textwidth]{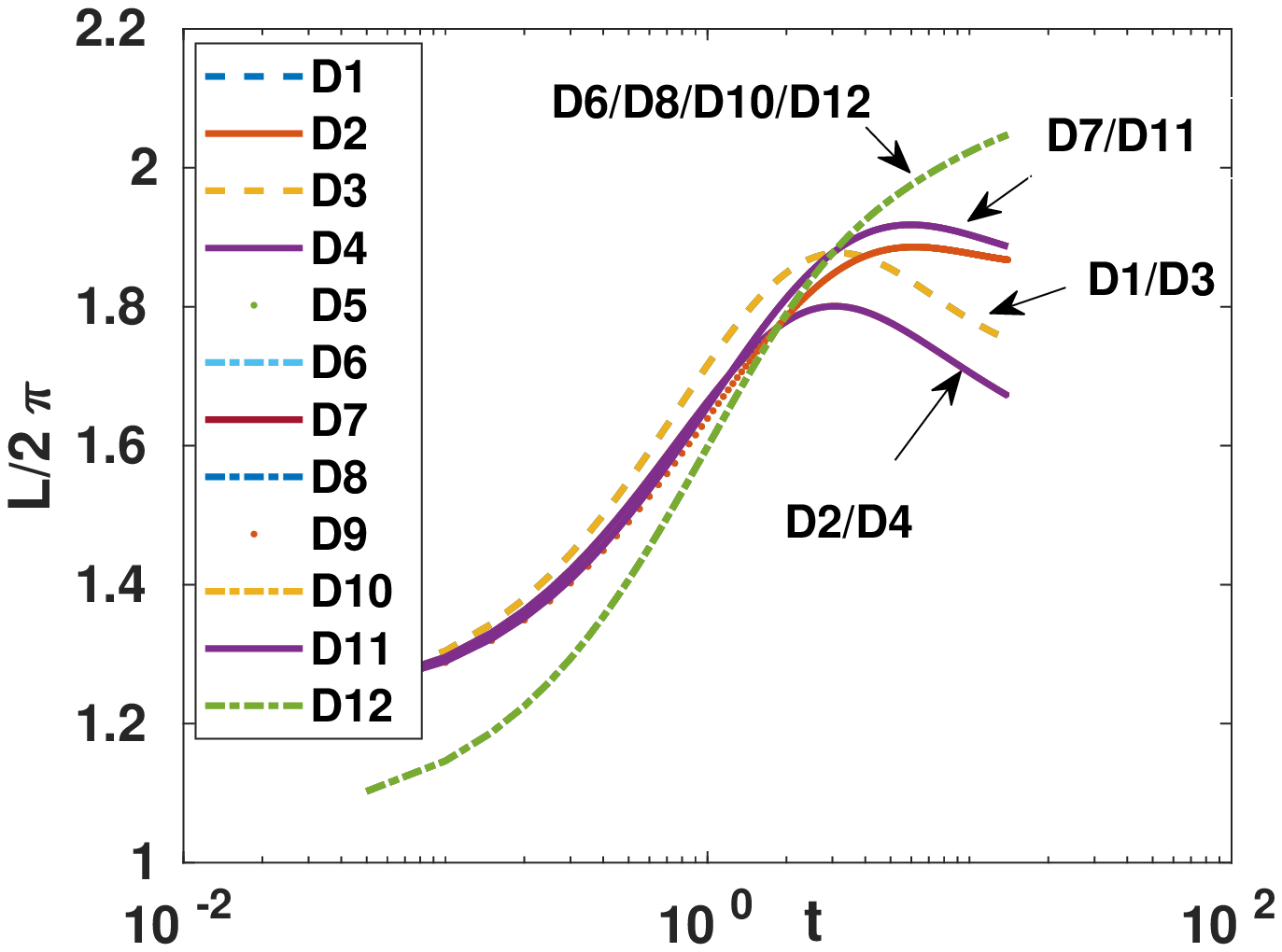}
\label{twelve_fig2D} }
\end{subfigure}
\caption{Time evolution of 12 elliptic regions D1 to D12 with centroids at $\left(3.75,0\right)$, 
$\left(0,4\right)$, 
$\left(-3.75,0\right)$, 
$\left(0,-4\right)$, 
$\left(7.5,0\right)$, 
$\left(5,5\right)$, 
$\left(0,-7\right)$, 
$\left(-5,5\right)$, 
$\left(-7.5,0\right)$, 
$\left(-5,-5\right)$, 
$\left(0,-7\right)$, 
and $\left(5,-5\right)$ at $t=0$, respectively. The domains D6, D8, D10, and D12 have major and minor axes $a=1.2$ and $b=0.9$ while the rest of the domains have major and minor axes $a=1.5$ and $b=0.9$.  We set $R_{\infty}=9$.}
\end{figure}

In the last simulation, we present the results of a simulation with a twelve domain configuration in Figs.~\ref{twelve_fig2B}, ~\ref{twelve_fig2C}, and ~\ref{twelve_fig2D}. The domains are arranged in ``two rings''. The inner ring consists of four particles (D1 to D4 arranged in counter clockwise direction having centroids at $\left(3.75,0\right)$, $\left(0,4\right)$, $\left(-3.75,0\right)$, and $\left(0,-4\right)$, respectively)
and the outer ring consists of eight particles, D5 to D12. Their centroids are located at $\left(7.5,0\right)$, $\left(5,5\right)$, $\left(0,-7\right)$, $\left(-5,5\right)$, $\left(-7.5,0\right)$, $\left(-5,-5\right)$, $\left(0,-7\right)$, and $\left(5,-5\right)$, respectively. The initial configuration has several symmetries which are all preserved in the simulation. The configuration enters the equilibrium phase at $t_c=11.3$ and does not show any coarsening type behaviour up to $t=14$. We observe that the domains in the outer ring grows more than the domains in the inner ring. This is probably due to the initial geometry where the outer domains have more space to grow and the inner domains are ``squeezed'' by the outer ring. Going by our previous simulation, we believe that placement of a central ellipse at $\left(0,0\right)$ will result in coarsening. 

\section{Summary and Conclusion}\label{concl}
In this article, we derived and studied a limiting case of Ohta-Kawasaki model. The resulting model is a variant of the Hele-Shaw problem. We then investigated the equations of the model using linear analysis and we reformulated the problem as boundary integral equations. Using small-scale decomposition technique for the equation of dynamics, we ran numerical simulations of these equations using a spectrally accurate algorithm in space and a second-order accurate temporal scheme. We investigated, with our numerical simulations, the evolution of different configurations of phase domains. Our simulations captured accurately the intermediate dynamics and final steady-state configuration, and reveals information about the far-field Dirichlet condition that drives the evolution. 

Choksi et al.~\cite{choksi2003derivation} related the Ohta-Kawasaki density functional theory (DFT) to the self-consistent mean field theory (SCFT) and \cite{mueller_rey_2018} compared the results of numerical simulations for the DFT, SCFT, and the Swift-Hohenberg model. Our future work will build upon these studies and the results introduced in this paper by comparing numerical simulations from the DFT, SCFT and the boundary integral method. Specifically, the energies of the stationary states and the metastability of the defect structures of the three models will be investigated. This will establish the feasibility of the boundary integral method for phase space exploration.

\clearpage
\begin{appendix}
\setcounter{equation}{0}
\renewcommand\theequation{A.\arabic{equation}}
\section{Derivation of the sharp-interface model } \label{sec:sil_derivation}
\subsection{Outer expansions} \label{subsec:sil_outerexpansion}
We assume $\phi(\tau,\vec{x})$, $\mu(\tau,\vec{x})$ and $\psi(\tau,\vec{x})$ have the asymptotic expansions,
$\phi  = \phi_0 + \varepsilon \phi_1 + \varepsilon^2 \phi_2 + 
\mathcal{O}(\varepsilon^3),
\mu  = \mu_0 + \varepsilon \mu_1 + \varepsilon^2 \mu_2 + 
\mathcal{O}(\varepsilon^3), 
\psi  = \psi_0 + \varepsilon \psi_1 + \varepsilon^2 \psi_2 + 
\mathcal{O}(\varepsilon^3)$.
The asymptotic problems in ``outer'' variables are for $\phi_i$ 
\begin{equation}
\mathcal{O}\left(\varepsilon^{0}\right): \partial_\tau \phi_0 = \Delta \mu_0,\qquad
\mathcal{O}\left(\varepsilon^{1}\right): \partial_\tau \phi_1 = \Delta \mu_1,\qquad
\mathcal{O}\left(\varepsilon^{2}\right): \partial_\tau \phi_2 = \Delta \mu_2.
\end{equation}
Similarly for $\mu_i$,
\begin{subequations}
\begin{alignat}{2}
&\mathcal{O}\left(\varepsilon^{-1}\right):
\qquad 0 &&= F^\prime \left( \phi_0 \right),
\label{eqn:sil_outerChem1_apx}\\
&\mathcal{O}\left(\varepsilon^{0}\right):
\qquad \mu_0 &&= F^{\prime\prime} \left( \phi_0 \right) \phi_1  + \psi_0,
\label{eqn:sil_outerChem2_apx}\\
&\mathcal{O}\left(\varepsilon^{1}\right):
\qquad \mu_1 &&= F^{\prime\prime} \left( \phi_0 \right) \phi_2 + \frac{1}{2} F^{\prime\prime\prime}\left(\phi_0\right)\phi_1^2 - \Delta \phi_0 + \psi_1.
\label{eqn:sil_outerChem3_apx}
\end{alignat}
\end{subequations}
and $\psi_i$,
\begin{equation}
\mathcal{O}\left(\varepsilon^{0}\right): - \Delta \psi_0 = \phi_0 - \bar{\phi},\quad
\mathcal{O}\left(\varepsilon^{1}\right): - \Delta \psi_1 = \phi_1,\quad
\mathcal{O}\left(\varepsilon^{2}\right): - \Delta \psi_2 = \phi_2.
\end{equation}
On the fixed boundary $\partial\Omega$, the boundary conditions for the asymptotic subproblems are
\begin{equation*}
\frac{\partial\phi_i}{\partial n_{\partial\Omega}} = 0, \qquad \frac{\partial \mu_i}{\partial n_{\partial\Omega}} = 0, \qquad \frac{\partial\psi_i}{\partial n_{\partial\Omega}} = 0, \quad \text{on} \quad \partial\Omega \quad \text{for}, \quad i=0,1,2,\dots
\end{equation*}

\subsection{Inner-outer coordinate transformations\label{subsec:sil_coordinatetransformaton}}
To derive the inner problems it is convenient to introduce a parametrization $\vec{r}(\tau,s) = (r_1(\tau,s),r_2(\tau,s))$ of the free interface, i.e. the sharp interface $\Gamma$ via the arc length $s$, and $\vec{\nu}(\tau,s)$ the normal inward-pointing vector along the free boundary, so that any point in the thin $\varepsilon$-region around $\Gamma$ can be expressed by  
\begin{equation*}
\vec{x}(\tau,s,z) = \vec{r}(\tau,s) + \varepsilon z \vec{\nu}(\tau,s).
\end{equation*}
where $\varepsilon z$ is the distance alongthe inward normal direction $\vec{\nu}(\tau,s)$ from the sharp interface $\Gamma$, given by 
%
%
\begin{equation*}
\vec{\nu}(\tau,s) = \left( -\partial_s r_2, \partial_s r_1 \right),\quad
\vec{t}(\tau,s) = \left( \partial_s r_1, \partial_s r_2 \right).
\end{equation*}
The relation the derivatives of a quantity $\tilde{v}(\tau,s,z)$ defined in inner coordinates to derivatives in the outer coordinates $v(\tau,\vec{x})$ can be expresses as a multiplication of matrices,
\begingroup
\begin{equation*}
\begin{bmatrix}
\partial_s \tilde{v}\\
\partial_z \tilde{v}\\
\partial_{\tau} \tilde{v}
\end{bmatrix}
= 
\begin{bmatrix}
\partial_s x & \partial_s y & 0 \\
\partial_z x & \partial_z y & 0 \\
\partial_{\tau} x & \partial_{\tau} y & 1
\end{bmatrix}
\cdot
\begin{bmatrix}
\partial_x v\\
\partial_y v\\
\partial_{\tau} v
\end{bmatrix},
\label{eqn:outerininner}
\end{equation*}
\endgroup
and vice versa
\begingroup
\renewcommand*{\arraystretch}{1.1}
\begin{equation*}
\begin{bmatrix}
\partial_x v\\
\partial_y v\\
\partial_{\tau} v
\end{bmatrix}
= 
\begin{bmatrix}
(1+\varepsilon z \kappa) \partial_s r_1 & 
- \varepsilon^{-1} \partial_s r_2 & 0 \\
(1+\varepsilon z \kappa) \partial_s r_2 & 
\varepsilon^{-1} \partial_s r_1 & 0 \\
-(1+\varepsilon z \kappa)V^{\vec{t}} & - \varepsilon^{-1} V^{\vec{\nu}} & 1
\end{bmatrix}
\cdot
\begin{bmatrix}
\partial_s \tilde{v}\\
\partial_z \tilde{v}\\
\partial_{\tau} \tilde{v}
\end{bmatrix},
\label{eqn:outer-to-inner}
\end{equation*}
\endgroup
where
\begin{equation*}
V^{\vec{t}} = \partial_{\tau} \vec{x} \cdot \vec{t}, \quad \text{and} \quad V^{\vec{\nu}} = \partial_{\tau} \vec{x} \cdot \vec{\nu},
\end{equation*}
denote the tangential and normal velocity of the free boundary respectively, with $\kappa$ denoting the curvature of the free boundary. Thus, the expression of the rescaled time derivative in terms of the inner-coordinates,
\begin{equation*}
\partial_\tau v = -(1+ \varepsilon z \kappa) \partial_s \tilde{v} - \varepsilon^{-1} V^{\vec{\nu}} \partial_z \tilde{v} + \partial_\tau \tilde{v}.
\end{equation*}
Applying the respective derivatives to higher order yields 
\begin{subequations}
\begin{align*}
\partial_{xx} v &= \varepsilon^{-2} \left( \partial_s r_2 \right)^2 \partial_{zz} \tilde{v} - \varepsilon^{-1} \left[ \kappa \left( \partial_s r_1 \right)^2 \partial_z \tilde{v} + 2 \partial_s r_1 \partial_s r_2 \partial_{sz} \tilde{v} \right] \nonumber\\
&
+ \left( \partial_s r_1 \right)^2 \partial_{ss} \tilde{v} - 2 \kappa \partial_s r_1 \partial_s r_2 \partial_s \tilde{v} 
- z \kappa \left[ \kappa \left( \partial_s r_1 \right)^2 \partial_z \tilde{v} + 2 \partial_s r_1 \partial_s r_2 \partial_{sz} \tilde{v} \right],\\
\partial_{yy} v &= \varepsilon^{-2} \left( \partial_s r_1 \right)^2 \partial_{zz} \tilde{v} - \varepsilon^{-1} \left[ \kappa \left( \partial_s r_2 \right)^2 \partial_z \tilde{v} - 2 \partial_s r_1 \partial_s r_2 \partial_{sz} \tilde{v} \right] \nonumber\\
&
+ \left( \partial_s r_2 \right)^2 \partial_{ss} \tilde{v} + 2 \kappa \partial_s r_1 \partial_s r_2 \partial_s \tilde{v} 
- z \kappa \left[ \kappa \left( \partial_s r_2 \right)^2 \partial_z \tilde{v} - 2 \partial_s r_1 \partial_s r_2 \partial_{sz} \tilde{v} \right]
\end{align*}
and for the Laplace operator in the inner-coordinates,
\begin{align*}
\Delta v = \partial_{xx} v + \partial_{yy} v 
= \varepsilon^{-2} \partial_{zz} \tilde{v} - \varepsilon^{-1} \kappa \partial_z \tilde{v} + \partial_{ss} \tilde{v} - z \kappa^2 \partial_z \tilde{v
}.
\end{align*}
\end{subequations}

\subsection{Inner expansions\label{subsec:sil_innerexpansion}}
We assume that inner asymptotic expansions for $\tilde{\phi}(\tau, s, z)$, $\tilde{\mu}(\tau, s,z )$ and $\tilde{\psi}(\tau, s,z )$ are given by
$\tilde{\phi} = \tilde{u}_0 + \varepsilon \tilde{u}_1 + \varepsilon^2 \tilde{u}_2 + \mathcal{O}(\varepsilon^3)$, 
$\tilde{\mu} = \tilde{\mu}_0 + \varepsilon \tilde{\mu}_1 + \varepsilon^2 \tilde{\mu}_2 + \mathcal{O}(\varepsilon^3)$, 
$\tilde{\psi} = \tilde{\psi}_0 + \varepsilon \tilde{\psi}_1 + \varepsilon^2 \tilde{\psi}_2 + \mathcal{O}(\varepsilon^3)$.
Application of the coordinate transformations to the governing equations yields the  asymptotic subproblems for the inner region 
for $\tilde\phi$ up till $\mathcal{O}(\varepsilon^0)$,
\begin{subequations}
\begin{alignat}{3}
&\mathcal{O}\left(\varepsilon^{-2}\right): 
\quad &0 &= \partial_z^2 \tilde{\mu}_0,
\label{eqn:sil_innerDiff1}\\
&\mathcal{O}\left(\varepsilon^{-1}\right):
\quad &-V^{\vec{\nu}} \partial_z \tilde{\phi}_0 &=  \partial_z^2 \tilde{\mu}_1 - \kappa \partial_z \tilde{\mu}_0,
\label{eqn:sil_innerDiff2}\\
&\mathcal{O}\left(\varepsilon^{0}\right):
\quad &-\partial_z \tilde{\phi}_0 - V^{\vec{\nu}} \partial_z \tilde{\phi}_1 + \partial_\tau \tilde{\phi}_0 &= \partial_z^2 \tilde{\mu}_2 - \kappa \partial_z \tilde{\mu}_1 + \partial_s^2 \tilde{\mu}_0 - z \kappa^2 \partial_z \tilde{\mu}_0.
\label{eqn:sil_innerDiff3}
\end{alignat}
\end{subequations}
For the chemical potential $\tilde\mu$ up to $\mathcal{O}(\varepsilon)$,
\begin{subequations}
\begin{alignat}{2}
&\mathcal{O}\left(\varepsilon^{-1}\right): \quad
0 &&= F^\prime(\tilde{u}_0) - \partial_z^2 \tilde{u}_0,\hspace{2.0cm}
\label{eqn:sil_innerChem1}\\
&\mathcal{O}\left(\varepsilon^{0}\right): \quad
\tilde{\mu}_0 &&= F^{\prime\prime}(\tilde{u}_0)\tilde{u}_1 + \kappa \partial_z \tilde{u}_0 - \partial_z^2 \tilde{u}_1 + \tilde{w}_0,
\label{eqn:sil_innerChem2}\\
&\mathcal{O}\left(\varepsilon^{1}\right): \quad
\tilde{\mu}_1 &&= - \partial_z^2 \tilde{u}_2 + \kappa \partial_z \tilde{u}_1 - \partial_s^2 \tilde{u}_0 + z \kappa^2 \partial_z \tilde{u}_0 + F^{\prime\prime}(\tilde{u}_0)\tilde{u}_2 + \frac{1}{2}F^{\prime\prime\prime}(\tilde{u}_0)\tilde{u}_1^2 \nonumber \\
& && \qquad + \tilde{w}_1,
\label{eqn:sil_innerChem3}
\end{alignat}
\end{subequations}
and for $\tilde\psi$,
\begin{subequations}
\begin{alignat}{2}
&\mathcal{O}\left(\varepsilon^{-2}\right): 
\qquad &-\partial_z^2 \tilde{\psi}_0 &= 0,
\label{eqn:sil_innerNL1}\\
&\mathcal{O}\left(\varepsilon^{-1}\right): 
\qquad &-\partial_z^2 \tilde{\psi}_1 + \kappa \partial_z \tilde{\psi}_0 &= 0,
\label{eqn:sil_innerNL2}\\
&\mathcal{O}\left(\varepsilon^{0}\right): 
\qquad &-\partial_z^2 \tilde{\psi}_2 + \kappa \partial_z \tilde{\psi}_1 - \partial_s^2 \tilde{\psi}_0 + z \kappa \partial_z \tilde{\psi}_0 &= \tilde{\phi}_0 - \bar{\phi}.
\label{eqn:sil_innerNL3}
\end{alignat}
\end{subequations}

\subsection{Matching\label{subsec:sil_results}}

From the leading order problem of the inner expansion for the chemical potential subequation (\ref{eqn:sil_innerChem1}),
\begin{equation*}
F^\prime(\tilde{\phi}_0) - \partial_z^2 \tilde{\phi}_0 = 0.
\end{equation*}
Multiplying by $\partial_z \tilde{\phi}_0$ and integrating in $z$ from $-\infty$ to $\infty$, 
\begin{equation*}
\int_{\phi_0^-}^{\phi_0^+} F^\prime(\tilde{\phi}_0) \dd \tilde{\phi}_0 = \int_{-\infty}^{\infty} \left(\partial_z^2 \tilde{\phi}_0\right)\partial_z \tilde{\phi}_0 \dd z,
\end{equation*}
where the integration limits are $\lim_{z \rightarrow \pm \infty} \tilde{\phi}_0(\tau,s,z) = \phi_0^\pm$ respectively. 
Since $\lim_{z \rightarrow \pm \infty} \frac{\partial \tilde{\phi}_0}{\partial z} =0$ for $\tilde{\phi}_0$ to be bounded. This leaves
\begin{equation*}
\int_{\phi_0^-}^{\phi_0^+} F^\prime \left( \tilde{\phi}_0 \right) \dd \tilde{\phi}_0  = 0,
\end{equation*}
which states that for the symmetric double-well potential 
the $x$-axis corresponding to $F^\prime \left( \tilde{\phi}_0 \right) = 0$ is the line of intersection that divides $F^\prime(\tilde{\phi}_0)$ such that the areas below and above the curve are equal. This implies that the limits of the integral are the points of intersection, i.e.
\begin{equation}
\phi_0^\pm = \pm 1 \quad \text{in} \quad {\Omega}^{\pm} \quad \text{resp.}
\label{eqn:sil_result1}
\end{equation}
This implies for the leading order outer problem in $\mu$
\begin{equation}
\Delta \mu_0 = 0 \qquad \text{in} \qquad \Omega \backslash \Gamma 
\label{eqn:sil_result2}
\end{equation}
and for $\psi_0$
\begin{equation}
\Delta\psi_0  = - (\phi_0 - \bar{\phi}),
\label{eqn:sil_result3}
\end{equation}
with boundary conditions
\begin{equation}
\frac{\partial \phi_0}{\partial n_{\partial\Omega}} = 0, \qquad \frac{\partial \mu_0}{\partial n_{\partial\Omega}} = 0, \qquad \frac{\partial \psi_0}{\partial n_{\partial\Omega}} = 0, \quad \text{on} \quad \partial {\Omega}.
\label{eqn:sil_result4}
\end{equation}


To proceed with the matching we write down the matching conditions by expanding inner and outer expansions, and express one of them (here the outer) in terms of the inner independent variables. Then we regroup in orders of $\varepsilon$ and obtain 
\begin{subequations}%
\begin{alignat}{2}
\mu_0^\pm &= \lim_{z \rightarrow \pm \infty} \tilde{\mu}_0 (\tau, \vec{r}, z),\label{eqn:sil_tildemu0matching}\\
\mu_1^\pm + z \vec{\nu} \cdot \nabla \mu^\pm_0   &= \lim_{z \rightarrow \pm \infty} \tilde{\mu}_1 (\tau, \vec{r}, z),\label{eqn:sil_tildemu1matching}\\
\mu^\pm_2 + z \vec{\nu} \cdot \nabla \mu_1^\pm + \frac{1}{2} z^2 \vec{\nu} \cdot \Delta \mu_0^\pm \cdot \vec{\nu}^\intercal \mu_0^\pm  &= \lim_{z \rightarrow \pm \infty} \tilde{\mu}_2 (\tau, \vec{r}, z).
\end{alignat}
\end{subequations}

Integrating (\ref{eqn:sil_innerDiff1}) twice gives
\begin{equation*}
\tilde{\mu}_0 = a_0 z + b_0.
\end{equation*}
Matching $\tilde{\mu}_0$ to $\mu_0^\pm$ by means of (\ref{eqn:sil_tildemu0matching}) yields $a_0 = 0$ and $\tilde{\mu}_0 = b_0=$constant. Next, notice that differentiating (\ref{eqn:sil_innerChem1}) with respect to $z$ and multiplying by $\tilde{\phi}_1$ yields
\begin{equation}
F^{\prime\prime}(\tilde{\phi}_0) \left(\partial_z \tilde{\phi}_0\right) \tilde{\phi}_1 - \left( \partial_z^3 \tilde{\phi}_0 \right)\tilde{\phi}_1 = 0.
\label{eqn:sil_fppRelations}
\end{equation}
Multiplying the next-order problem of the inner chemical potential  (\ref{eqn:sil_innerChem2}) by $\partial_z \tilde{\phi}_0$ and using (\ref{eqn:sil_fppRelations}) gives 
\begin{equation*}
\tilde{\mu}_0 \left(\partial_z \tilde{\phi}_0\right) = \left( \partial_z^3 \tilde{\phi}_0 \right)\tilde{\phi}_1 + \kappa \left(\partial_z \tilde{\phi}_0\right)^2 - \left(\partial_z^2 \tilde{\phi}_1\right)\left(\partial_z \tilde{\phi}_0\right) + \tilde{\psi}_0 \left(\partial_z \tilde{\phi}_0\right).
\end{equation*}
Integrating the above in $z$ from $-\infty$ to $\infty$, applying integration by parts and using the boundedness of the leading order $\tilde{\phi}_0$ and the leading order non-local term $\tilde{\psi}_0$ is a functional of $\tilde{\phi}_0$ with $\lim_{z \rightarrow \pm \infty} \tilde{\psi}_0 = \psi_0[\phi_0^\pm]$ we obtain 
\begin{equation*}
\tilde{\mu}_0 \left[ \tilde{\phi}_0 \right]_{-\infty}^{\infty} = \kappa \int_{-\infty}^{\infty} \left(\partial_z \tilde{\phi}_0\right)^2 \dd z  + \tilde{\psi}_0 \left[ \tilde{\phi}_0 \right]_{-\infty}^{\infty},
\end{equation*}
where $\int_{-\infty}^{\infty} \partial_z \tilde{\phi}_0 \dd z = \left[ \tilde{\phi}_0 \right]_{-\infty}^{\infty}$, the jump of $\tilde{\phi}_0$ over the interface. Dividing by $\left[ \tilde{\phi}_0 \right]_{-\infty}^{\infty}$  and setting
%
\begin{equation*}
\frac{\int_{-\infty}^{\infty} \left(\partial_z \tilde{\phi}_0\right)^2 \dd z}{\left[ \tilde{\phi}_0 \right]_{-\infty}^{\infty}} = C,
\end{equation*}
which is a constant, we obtain 
\begin{equation*}
\tilde{\mu}_0 = C \kappa + \tilde{\psi}_0.
\end{equation*}
The next-order matching conditions then implies 
%
\begin{equation}
\mu_0 = C \kappa + \psi_0 \qquad \text{on} \qquad \Gamma.
\label{eqn:sil_result5}
\end{equation}


To obtain the normal velocity of the free boundary $V^{\vec{\nu}}$
we integrate (\ref{eqn:sil_innerDiff2}) from $-\infty$ to $\infty$, 
\begin{equation}
-V^{\vec{\nu}} =  \frac{1}{2} \underbrace{\left[ \partial_z \tilde{\mu}_1 \right]_{-\infty}^{\infty}}_{\raisebox{.5pt}{\textcircled{\raisebox{-.9pt} {A}}}} 
- \frac{1}{2} \kappa \underbrace{\left[ \tilde{\mu}_0 \right]_{-\infty}^{\infty}}_{\raisebox{.5pt}{\textcircled{\raisebox{-1.2pt} {B}}}},
\label{eqn:sil_vIntermediate}
\end{equation}
From \ref{eqn:sil_result5}, $\tilde{\mu}_0$ is independent of $z$, so ${\raisebox{.5pt}{\textcircled{\raisebox{-1.2pt} {B}}}} = 0$. Furthermore, notice that differentiating the next-order matching of $\tilde{\mu}_1$ in (\ref{eqn:sil_tildemu1matching}) with respect to $z$ yields
\begin{align*}
\left. \partial_z \tilde{\mu}_1 \right\rvert_{z=-\infty}^{\infty} &= \underbrace{\left. \partial_z \mu_1 \right\rvert_-^+}_{=0} + \left. \vec{\nu} \cdot \nabla \mu_0 \right\rvert_-^+ \left. \vec{\nu} \cdot \nabla \mu_0 \right\rvert_-^+ \equiv {\raisebox{.5pt}{\textcircled{\raisebox{-.9pt} {A}}}},
\end{align*}
with $\left. \partial_z \mu_1 \right\rvert_-^+ = 0$ since the outer $\mu^\pm_1$'s are independent of $z$. Substituting these results back into (\ref{eqn:sil_vIntermediate}),
\begin{align}
V^{\vec{\nu}} &= - \frac{1}{2} \left[ \frac{\partial \mu_0}{\partial \vec{\nu}} \right]_\Gamma.
\label{eqn:sil_result6}
\end{align}
\end{appendix}
\section*{Acknowledgments}
R. C. thanks the Deutsche Forschungsgemeinschaft (DFG) for the funding through CRC 1114 ``Scaling Cascades in Complex Systems'' (Project Number 235221301), Project A02 and the Weierstrass Institute. S. L. and J. L. gratefully acknowledge partial support from the National Science Foundation, Division of Mathematical Sciences through grants NSF-DMS 1719960 (J. L.) and NSF-DMS 1720420 (S. L.).  S. L. is also partially supported by grant  NSF-ECCS 1927432. J.L. also acknowledges partial support from grants NSF-DMS 1763272 and the Simons Foundation (594598, QN) for the Center for Multiscale Cell Fate Research at UC Irvine.
\bibliographystyle{elsarticle-num}
\bibliography{jcp_format}
\end{document}